\documentclass[11pt,a4paper,bibtotoc,final]{scrartcl}
\usepackage{a4wide}
\usepackage{amsfonts}
\usepackage{amsmath}
\usepackage{amssymb}
\usepackage{amsthm}
\usepackage[section]{algorithm}
\usepackage{algorithmic}
\usepackage{pgfplots}
\usepackage{graphicx}
\usepackage{color}
\usepackage{colortbl}
\usepackage{scrpage2}
\usepackage{multirow}
\usepackage{boxedminipage}
\usepackage{booktabs}
\usepackage{enumerate}
\usepackage[caption=false,font=footnotesize]{subfig}

\usepackage{todonotes}
\usepackage{placeins}

\newcommand{\N}{\ensuremath{\mathbb{N}}}

\newcommand{\T}{\ensuremath{\mathbb{T}}}

\newcommand{\Z}{\ensuremath{\mathbb{Z}}}

\newcommand{\R}{\ensuremath{\mathbb{R}}}
\newcommand{\C}{\ensuremath{\mathbb{C}}}

\newcommand{\ii}{\textnormal{i}}
\newcommand{\e}{\textnormal{e}}

\newcommand{\ceil}[1]{\left\lceil#1\right\rceil}
\newcommand{\floor}[1]{\left\lfloor#1\right\rfloor}
\newcommand{\zb}[1]{\ensuremath{\boldsymbol{#1}}}

\newcommand{\AC}{\ensuremath{\operatorname{X}}}
\newcommand{\vol}{\ensuremath{\operatorname{vol}}}

\newcommand{\bigtimes}{\mathop{\text{\Large{$\times$}}}}

\newcommand{\boldk}{{\ensuremath{\boldsymbol{k}}}}
\newcommand{\boldj}{{\ensuremath{\boldsymbol{j}}}}

\newcommand{\boldh}{{\ensuremath{\boldsymbol{h}}}}

\newcommand{\boldx}{{\ensuremath{\boldsymbol{x}}}}

\newcommand{\boldz}{{\ensuremath{\boldsymbol{z}}}}

\newcommand{\boldzero}{{\ensuremath{\boldsymbol{0}}}}

\newtheorem{theorem}{Theorem}[section]
\newtheorem{lemma}[theorem]{Lemma}
\newtheorem{remark}[theorem]{Remark}
\newtheorem{generalisation}[theorem]{Generalisation}
\newtheorem{definition}[theorem]{Definition}
\newtheorem{example}[theorem]{Example}
\newtheorem{corollary}[theorem]{Corollary}
\newtheorem{proposition}[theorem]{Proposition}

\newenvironment{Theorem}{\goodbreak \begin{theorem}\sl}{\end{theorem}}
\newenvironment{Lemma}{\goodbreak \begin{lemma}\sl}{\end{lemma}}
\newenvironment{Remark}{\goodbreak \begin{remark}\rm}{\end{remark}}

\newenvironment{Proposition}{\goodbreak \begin{proposition}\rm}{\end{proposition}}

\makeatletter
\def\imod#1{\allowbreak\mkern10mu({\operator@font mod}\,\,#1)}
\makeatother

\numberwithin{equation}{section}
\numberwithin{table}{section}
\numberwithin{figure}{section}

\newcommand{\bend}{\hspace*{0ex} \hfill \hbox{\vrule height
    1.5ex\vbox{\hrule width 1.4ex \vskip 1.4ex\hrule  width 1.4ex}\vrule
    height 1.5ex}}

\long\def\symbolfootnote[#1]#2{\begingroup%
\def\thefootnote{\fnsymbol{footnote}}\footnote[#1]{#2}\endgroup}
\newcommand{\sspan}{\textnormal{span}}

\renewcommand{\mathbf}[1]{\ensuremath{\boldsymbol{#1}}}
\renewcommand{\textbf}[1]{{\ensuremath{\boldsymbol{#1}}}}
\renewcommand{\thefootnote}{\fnsymbol{footnote}}
\title{Tight error bounds for rank-1 lattice sampling in spaces of hybrid mixed smoothness}

\date{}
\author{
Glenn Byrenheid \footnotemark[2] \and
Lutz K\"ammerer\footnotemark[1] \and
Tino Ullrich\footnotemark[2] \and 
Toni Volkmer\footnotemark[1]}

\newif\ifshowextendedpaperversion
\showextendedpaperversiontrue

\begin{document}

\maketitle

\begin{abstract} We consider the approximate recovery of multivariate periodic functions from a discrete set of
function values taken on a rank-$s$ integration lattice.
The main result is the fact that any (non-)linear
reconstruction algorithm taking function values on a rank-$s$ lattice of size $M$ has a dimension-independent
lower bound of $2^{-(\alpha+1)/2} M^{-\alpha/2}$ when considering the optimal worst-case error with respect to function spaces of
(hybrid) mixed smoothness $\alpha>0$ on the $d$-torus.
We complement this lower bound 
with upper bounds that coincide up to logarithmic terms.
These upper bounds are obtained by a detailed analysis of
a rank-1 lattice sampling strategy, where the rank-1 lattices
are constructed by a component--by--component (CBC) method.
This improves on earlier results obtained in \cite{KuSlWo06} and \cite{KuWaWo09}.
The lattice (group) structure allows for an 
efficient approximation of the underlying function from its sampled values
using a single one-dimensional fast Fourier transform.
This is one reason why these algorithms keep attracting significant interest.
We compare our results to recent (almost) optimal methods based upon samples on sparse grids.

\medskip


\medskip


\end{abstract}
\footnotetext[1]{
  Technische Universit\"at Chemnitz, Faculty of Mathematics, 09107 Chemnitz, Germany\\
  \{lutz.kaemmerer, toni.volkmer\}@mathematik.tu-chemnitz.de
}
\footnotetext[2]{
  Hausdorff Center for Mathematics, Endenicher Allee 62, 53115 Bonn, Germany\\
  \{glenn.byrenheid, tino.ullrich\}@hcm.uni-bonn.de
}

\medskip

\ifshowextendedpaperversion
\newpage
\fi
\newpage

\section{Introduction}
This paper deals with the reconstruction of multivariate periodic functions from a discrete set of $M$ function values
along rank-$1$ lattices. Such lattices are widely used for the efficient numerical integration of multivariate periodic
functions since the 1950ies \cite{Ba59, Ko59, Nie78, SlJo94, CoNu07} and represent a well-distributed set of points in
$[0,1)^d$. A rank-$1$ lattice with $M\in\N$ points and
generating vector $\zb z \in\Z^d$ is given by 
 \begin{equation*}
 \Lambda(\zb z, M):=\left\{\frac{j}{M}\zb z\bmod \zb 1\colon
j=0,\ldots,M-1\right\}\,.
 \end{equation*}
In this paper we will show that restricting the set of available discrete information to samples from a rank-$s$
lattice, cf. \cite{SlJo94}, seriously affects the rate of convergence of a corresponding worst-case error with respect to classes of
functions with (hybrid) mixed smoothness $\alpha>0$. To be more precise, for any (possibly nonlinear) reconstruction
procedure from sampled values along rank-$s$ lattices we can find a function in the
periodic Sobolev spaces $\mathcal{H}^{\alpha}_{\mathrm{mix}}$
such that the $L_2(\T^d)$ mean square error is at least $2^{-(\alpha+1)/2} M^{-\alpha/2}$.
In contrast to that, it has been proved recently, that the
sampling recovery from (energy) sparse grids leads to much better convergence rates, namely $M^{-\alpha}$ in the main
term, see \cite{TW07, T08, BDSU14}.

Subsequently, we study particular reconstructing algorithms, which are based on the naive
approach of approximating the potentially ``largest'' Fourier coefficients (integrals) with the same rank-1 lattice
rule.
Despite the lacking asymptotical optimality, recovery from so-called reconstructing rank-1 lattices, cf. \cite{Kae2013, KaPoVo13}, has some striking advantages.

First, the matrix of the underlying linear system of equations has orthogonal columns
due to the group structure \cite{Be13} and the reconstructing property of the used rank-1 lattices.
Consequently, the computation is stable, cf. \cite{KaKuPo10, Kae2013}.

Second, the CBC strategy \cite[Tab.\ 3.1]{kaemmererdiss} provides a search method for
a reconstructing rank-1 lattice
which allows for the computation of the approximate Fourier coefficients 
belonging to frequencies lying on potentially unstructured sets.
Besides a basic structure, e.g. generalized hyperbolic crosses,
additional sparsity in the structure of the set of basis functions can be
easily incorporated and may considerably reduce the number of required samples, e.g. see \cite[Example 6.1]{KaPoVo13}.

Last, the approximate reconstruction can be efficiently performed 
using the sampled values of the underlying function and applying a single
one-dimensional fast Fourier transform, cf. Algorithm~\ref{alg:r1l_reconstruction_fct} and \cite{LiHi03, Be13}. This idea has already been investigated by many
authors including two of the present ones, see \cite{Tem86,KuSlWo06,KuSlWo08,KuWaWo09,KaPoVo13}.
The arithmetic complexity is $\mathcal{O}(M\log M)$, and thus almost linear in the number of used sampling values.

The above mentioned advantages
motivate
a refined error analysis for the upper bounds which results in the
observation that for the rank-$1$ lattice sampling the lower bound $M^{-\alpha/2}$ is
sharp in the main order. It is important to mention that the rate $M^{-\alpha/2}$ is present in any dimension $d\ge 2$.
Hence, the proposed naive but fast reconstruction algorithm is already more accurate than a comparable full tensor
grid in case $d>2$ yielding the order $M^{-\alpha/d}$. Moreover, the comparison to the mentioned sparse grid techniques
is not completely hopeless since neither the asymptotical behavior of the approximation error tells anything about small values of $M$ (so-called
preasymptotics), which is indeed relevant for practical issues, nor is the computational cost for computing the sparse
grid approximant completely reflected in the (optimal) main rate $M^{-\alpha}$, cf. \cite{KuSiUl13,KSU13,KMU15}. This is the
reason why these algorithms keep attracting more and more interest recently.
 
We consider the rate of convergence in the number of lattice points
$M$ of the worst-case error with respect to periodic Sobolev spaces with bounded
mixed derivatives in $L_2$. These classes are given by 
\begin{equation}\label{f1}
\mathcal{H}^{\alpha}_{\mathrm{mix}}(\T^d)=\Big\{f\in L_2(\T^d):
\|f|\mathcal{H}^{\alpha}_{\mathrm{mix}}(\T^d)\|^2:=\sum_{\|\zb m\|_{\infty}\leq
\alpha}\|D^{\zb m}f\|_2^2<\infty\Big\}\,,
\end{equation}
where $\alpha\in \N$ denotes the mixed smoothness of the space. In order to quantitatively assess the quality of the
proposed approximation,
we introduce specifically tailored minimal worst-case errors $g^{\mathrm{latt}_1}_M({\mathcal{F}},Y)$ with respect to the
function class $\mathcal{F}$ and the error in the norm of the function class $Y$. Our main result in case ${\mathcal{F}} =
\mathcal{H}^{\alpha}_{\mathrm{mix}}(\T^d)$ and $Y=L_2(\T^d)$ reads as
follows
\begin{equation*}
 M^{-\alpha/2} \lesssim g^{\mathrm{latt}_1}_M(\mathcal{H}^{\alpha}_{\mathrm{mix}}(\T^d),L_2(\T^d))
\lesssim M^{-\alpha/2}(\log M)^{\frac{d-2}{2}\alpha+\frac{d-1}{2}}\,,\quad M\in \N\,.
\end{equation*}
To be more precise, we use the following definition for sampling numbers along 
rank-$1$ lattice nodes
\begin{align*}
g^{\mathrm{latt}_1}_M({\mathcal{F}},Y):=\inf_{\zb
z\in\Z^d}\operatorname{Samp}_{\Lambda(\zb 
z,M)}({\mathcal{F}},Y)\,,\quad M\in \N,
\end{align*}
where we put for $\mathcal{G} := \{\mathbf{x}^1,...,\mathbf{x}^M\}\subset \T^d$
\begin{align*}
\operatorname{Samp}_{\mathcal{G}}({\mathcal{F}},Y):=\inf_{
A:\C^M\to Y}\sup_{\|f|{\mathcal{F}}\|\leq 1}
\left\|f-A\big(f(\mathbf{x}^i)\big)_{i=1}^M\right\|_Y\,.
\end {align*}
Here we allow as well non-linear reconstruction operators $A\colon\C^M\to
Y$. The general (non-linear) sampling numbers are defined as
\begin{align*}
g_M({\mathcal{F}},Y):=\inf_{\mathcal{G}}\operatorname{Samp}_
{\mathcal{G}}({\mathcal{F}},Y)\,,\quad M\in \N,
\end{align*}
for arbitrary sets of sampling nodes $\mathcal{G} := \{\mathbf{x}^1,...,\mathbf{x}^M\}\subset \T^d$ and are sometimes
also referred to as ``optimal sampling recovery''. These quantities are not the central focus of this paper, they
rather serve as benchmark quantity. If the reconstruction operator $A$ is supposed to be linear then we will use the
notation $g_M^{\mathrm{lin}}({\mathcal{F}},Y)$. These quantities are well studied up to
some prominent logarithmic gaps (cf. 3rd column in Table \ref{tab:mainres:overview_general_r1l},
\ref{tab:mainres:overview_general_r1l_2} and \ref{tab:mainres:overview_fibonacci_r1l}). For an overview we refer to
\cite{BDSU14} and the references therein. 
Additionally, let us mention the work of Temlyakov \cite{Tem93_2,Tem93}, Griebel et
al. \cite{Bungartz.Griebel:2004,Griebel.Hamaekers:2014, GrKn09}, Dinh
\cite{DD11,DD15,BDSU14} , Sickel \cite{S03,S06,TW07,SiUl09,BDSU14}, Ullrich
\cite{TW07,T08,SiUl09,BDSU14}.

The main goal of this paper is to study the quantities $g^{\mathrm{latt}_1}_M({\mathcal{F}},Y)$ in several
different approximation settings. At first, we measure the error in $Y = L_q(\T^d)$ with $2\leq q \leq
\infty$. In addition, we consider worst-case errors measured in isotropic Sobolev spaces
$Y = \mathcal{H}^{\gamma}(\T^d)$ (defined as $\mathcal{H}^{\gamma}(\T^d):=\mathcal{H}^{0, \gamma}(\T^d)$ in
\eqref{eq:H_alpha_beta_integer} below) which includes the energy-norm $\mathcal{H}^{1}(\T^d)$ relevant for
Galerkin approximation schemes.
Multivariate functions are taken from fractional ($\alpha>0$) Sobolev
spaces ${\mathcal{F}} = \mathcal{H}^{\alpha}_{\mathrm{mix}}(\T^d)$ of mixed smoothness and even
more general hybrid type Sobolev spaces ${\mathcal{F}} =\mathcal{H}^{\alpha,\beta}(\T^d)$,
introduced by Griebel and Knapek \cite{GrKn09}. In fact, Yserentant \cite{Y00}
proved  that eigenfunctions of the positive spectrum of the electronic
Schr\"odinger
operators have a mixed type regularity. Even more, their regularity can be
described as a combination of mixed and isotropic (hybrid) smoothness
\begin{equation}\label{eq:H_alpha_beta_integer}
\mathcal{H}^{\alpha,\beta}(\T^d)=\Big\{f\in L_2(\T^d):
\|f|\mathcal{H}^{\alpha,\beta}(\T^d)\|^2:=\sum_{\|\zb m\|_{\infty}\leq
\alpha}\sum_{\|\zb n\|_{1}\leq \beta}\|D^{\zb m+\zb n}f\|_2^2<\infty\Big\}.
\end{equation}
A related concept is given by anisotropic mixed Sobolev smoothness
\begin{align}\label{anisomix}
\mathcal{H}^{\zb \alpha}_{\mathrm{mix}}(\T^d)=\Big\{f\in L_2(\T^d):
\|f|\mathcal{H}^{\zb \alpha}_{\mathrm{mix}}(\T^d)\|^2:=\sum_{\substack{m_i \leq \alpha_i\\i=1,\ldots,d}}\|D^{\zb
m}f\|_2^2<\infty\Big\}\,,
\end{align}
where the smoothness is characterized by vectors $\zb \alpha \in \N_0^d$. In fact, we have the representation
\begin{equation*}
    \mathcal{H}^{\alpha,\beta} = \bigcap\limits_{i=1}^d \mathcal{H}^{\alpha\cdot {\zb 1} + \beta\cdot
{\zb e}_i}_{\mathrm{mix}}\,,
\end{equation*}
where ${\zb e}_i$ is the $i$-th unit vector.
The norms in \eqref{f1}, \eqref{eq:H_alpha_beta_integer}, \eqref{anisomix} can be
rephrased as weighted $\ell_2$-sums of Fourier coefficients which is also
the natural way to extend the spaces $\mathcal{H}^{\alpha,\beta}(\T^d)$ to
fractional parameters, see \eqref{def:hab} below.
We extend methods from \cite{KaKuPo10,KaPoVo13} to obtain sharp bounds (up to logarithmic
factors) for $g^{\text{latt}_1}_M(\mathcal{H}^{\alpha,\beta}(\T^d),H^{\gamma}(\T^d))$,
which show in particular that even non-linear reconstruction maps can not get below
$c_{\alpha,\beta,\gamma,d}M^{-(\alpha+\beta-\gamma)/2}$. The upper bounds are obtained with a specific
simple algorithm that approximates the ``largest'' Fourier coefficients \eqref{eq:def_Fourier_coefficient} of the function
with one fixed lattice rule, where the corresponding frequencies of the Fourier coefficients are determined by the function class.
To this end, a so-called reconstructing rank-$1$ lattice \cite[Ch.\ 3]{kaemmererdiss} is used,
which is constructed via the component--by--component (CBC) strategy \cite{SlRe02}.
Similar strategies have already proved useful for numerical integration, see
\cite{SlRe02, CoKuNu10, CoNu07}. The basic idea behind is the construction of a generating vector $\zb
z$ component-wise by iteratively increasing the dimension of the index set for which a reproduction property should
hold.

Let us finally comment on some relevant earlier results in this direction. One of the first upper
bounds for $g^{\mathrm{latt}_1}_M(\mathcal{H}^{\alpha}_{\mathrm{mix}}(\T^d),L_2(\T^d))$
has been obtained by Temlyakov in \cite{Tem86} for the Korobov lattice, which
represents a rank-$1$
lattice with a generating vector $\mathbf{a} = (1,a,a^2,\ldots,a^{d-1})$ for
some integer $a$. He obtained the estimate
\begin{align*}
\operatorname{Samp}_{\Lambda(\mathbf{a},M)}(\mathcal{H}^{\alpha}_{\mathrm{mix}
}(\T^d),L_2(\T^d)) \lesssim M^{-\alpha/2} \,
(\log M)^{(d-1)(\alpha/2+1/2)}. 
\end{align*}
\begin{sloppypar}
\noindent Further results that imply upper bounds for
$g^{\mathrm{latt}_1}_M(\mathcal{H}^{\alpha}_{\mathrm{mix}}(\T^d),L_2(\T^d))$ have been proved in \cite{KuSlWo06}.
Rephrasing the error bounds in \cite{KuSlWo06} depending on the number of lattice points $M$, we observe a rate of
$M^{-(\alpha-\lambda)/2}$ for any $\lambda>0$.
In \cite{KuWaWo09} the rank-$1$ lattice sampling error measured in $L_{\infty}(\T^d)$ is considered and
the main rate $M^{-(\alpha-1/2-\lambda)/2}$ is obtained for every $\lambda>0$.
In \cite{KaPoVo14} the technique used by
Temlyakov \cite{Tem86} is expanded to model spaces $\mathcal{H}^{\alpha,\beta}(\T^d)$ with $\beta<0$ and
$\alpha+\beta>1/2$, where the authors obtain the upper bound
$$g^{\text{latt}_1}_M(\mathcal{H}^{\alpha,\beta}(\T^d),L_2(\T^d))\lesssim M^{-(\alpha+\beta)/2}$$
without any further logarithmic dependence.
\end{sloppypar}

{\bf Contribution and main results.}
The first main contribution of the present paper is the lower bound
$$
c_{\alpha,\beta,\gamma}\; M^{-(\alpha+\beta-\gamma)/2} \leq g^{\mathrm{latt}_1}_M(\mathcal{H}^{\alpha,\beta}(\T^d),Y),
\qquad c_{\alpha,\beta,\gamma} := 2^{-(\alpha+\beta-\gamma+1)/2},
$$
for
$Y\in\{L_2(\T^d)=\mathcal{H}^{0}(\T^d),\mathcal{H}^{\gamma}(\T^d),\mathcal{H}^{\gamma}_{\mathrm{mix}}(\T^d)\}$
and $\min\{\alpha,\alpha+\beta\}>\gamma\geq 0$, cf. Section~\ref{sec:lower_bounds}.
 In the cases $Y\in \{L_2(\T^d),\mathcal{H}^{\gamma}(\T^d),\mathcal{H}^{\gamma}_{\mathrm{mix}}(\T^d)\}$ and
$\alpha+\beta>\max\{\gamma,1/2\}$ with $\beta \le 0$ and $\gamma\ge0$, the
 upper bounds on the rank-$1$ lattice sampling rates
 match the general lower bounds up to logarithmic factors, cf. Sections \ref{sec:upper_bounds}, \ref{sec:2d}.
 \begin{table}[h]
 	\centering
 	\begin{tabular}{lll}
 		\toprule
 		Y & $g_M^{\mathrm{latt}_1}(\mathcal{H}^{\alpha,\beta}(\T^d),Y)$ &
 		$g_M(\mathcal{H}^{\alpha,\beta}(\T^d),Y)$ \\
 		\cmidrule{1-3}
 		$L_2(\T^d),\mathcal{H}^{\gamma}(\T^d)$,$\mathcal{H}^{\gamma}_{\mathrm{mix}}(\T^d)$&  $\gtrsim M^{-\frac{\alpha+\beta-\gamma}{2}}$ & $\gtrsim
 		M^{-(\alpha+\beta-\gamma)}$ \\
 		&{\footnotesize (Proposition \ref{prop:lower_bound_r1l_general})}& { \footnotesize\cite{DiUl12} linear,\;\cite{DD15} non-linear, non-periodic}
 		\\
 		\bottomrule
 	\end{tabular}
 	\caption{Lower bounds of sampling numbers for different sampling methods.}
 	\label{tab:mainres:lower_bounds}
 \end{table}
 
The second column in Table~\ref{tab:mainres:lower_bounds},
\ref{tab:mainres:overview_general_r1l} and
\ref{tab:mainres:overview_general_r1l_2} is headlined with 
$g_M^{\mathrm{latt}_1}(\mathcal{H}^{\alpha,\beta}(\T^d),Y)$ and  
presents lower and upper bounds on the sampling rates in various settings for sampling
along reconstructing rank-1 lattices.
Table~\ref{tab:mainres:lower_bounds} shows the lower bounds from Section~\ref{sec:lower_bounds}.
Table \ref{tab:mainres:overview_general_r1l} deals with upper bounds in the model spaces
$\mathcal{H}^{\alpha}_{\mathrm{mix}}(\T^d)$, whereas in Table \ref{tab:mainres:overview_general_r1l_2} model spaces
$\mathcal{H}^{\alpha,\beta}(\T^d)$ with negative isotropic smoothness parameter $\beta$ are considered. The
corresponding $L_2(\T^d)$ error estimate in the first table improves on the result obtained by Temlyakov in
\cite{Tem86}
by a logarithmic factor  $(\log M)^{\alpha/2}$. In contrast to the rank-$1$
lattices constructed by the CBC strategy, the considerations by Temlyakov are
based on rank-$1$ lattices of Korobov type. Smoothness parameters are chosen from
$\beta<0$, $\alpha+\beta>\max\{\gamma,1/2\}$, $\gamma>0$, and $2<q<\infty$. Best known bounds are based on energy
sparse grid sampling.
References marked with $^*$ mean that the result is not stated there explicitly but follows with the same method
therein.
For our method the crucial property of the used rank-1 lattice sampling scheme is the reconstruction property
\eqref{eqn:reco_prop1}.
In order to construct such rank-1 lattices, one may use the CBC strategy \cite[Tab.\ 3.1]{kaemmererdiss}. Additionally,
in case $d=2$ the Fibonacci
lattice fulfills such a property. In both of these cases, we obtain the improved estimates as shown in
Table~\ref{tab:mainres:overview_fibonacci_r1l}.
Smoothness parameters are chosen from
$\alpha>1/2$, $\alpha>\gamma > 0$. The upper bounds for $g_M^{\mathrm{latt}_1}$ are realized either by the
Fibonacci or CBC-generated lattice. From the point of error estimates, the case $d=2$ represents an interesting special
case. We have sharp bounds and no logarithmic
dependencies here, except in the case where we measure the error in a space with mixed regularity. Hence, lattice
sampling turns out to be as good as sampling on the full tensor grid in $d=2$. Last but not least, we consider the
recovery of functions from $\mathcal{H}^{\zb \alpha}_{\mathrm{mix}}(\T^d)$ with anisotropic mixed
smoothness. We consider
smoothness vectors $\zb \alpha\in\R^d$ with  first $\mu$ smallest smoothness directions, i.e.
$$\frac{1}{2}<\alpha_1=\ldots=\alpha_{\mu}<\alpha_{\mu+1}\leq \ldots\leq \alpha_d.$$ Here we show for the $L_{\infty}$
approximation error the bound
$$g^{\mathrm{latt}_1}_M(\mathcal{H}^{\zb \alpha}_{\mathrm{mix}}(\T^d),L_{\infty}(\T^d))\lesssim
M^{-(\alpha_1-\frac{1}{2})/2}
(\log M)^{(\mu-1)(\alpha_1/2+1/4)}.$$
That means the exponent of the logarithm depends only on $\mu<d$ instead of $d$. Similar effects are also known for
general linear approximation and sparse grid sampling, cf. \cite{DD15,Tem87}.
\begin{table}[h]
	\begin{center}
		\begin{tabular}{lll}
			\toprule
			$Y$ & $g_M^{\mathrm{latt}_1}(\mathcal{H}^{\alpha}_{\text{mix}}(\T^d),Y)$ &
			$g_M^{\mathrm{lin}}(\mathcal{H}^{\alpha}_{\text{mix}}(\T^d),Y)$ \\
			\cmidrule{1-3}
			$L_2(\T^d)$  & $\lesssim M^{-\frac{\alpha}{2}}(\log
			M)^{\frac{d-2}{2}\alpha+\frac{d-1}{2}}$ & $\lesssim M^{-\alpha}(\log
			M)^{(d-1)(\alpha+\frac{1}{2})}$ \\
			&{\footnotesize (Theorem \ref{thm:alias_err:dyadic})}&{
				{\footnotesize\cite[Theorem 6.10]{BDSU14}, sparse grid}}\\
			\cmidrule{1-3}
			$L_q(\T^d)$  & $\lesssim M^{-\frac{\alpha-(\frac{1}{2}-\frac{1}{q})}{2}}$ & $\asymp
			M^{-(\alpha-(\frac{1}{2}-\frac{1}{q}))}$ \\
			&$\quad(\log	M)^{\frac{d-2}{2}(\alpha-(\frac{1}{2}-\frac{1}{q}))+\frac{d-1}{2}}$
			&$\quad(\log	M)^{(d-1)(\alpha-(\frac{1}{2}-\frac{1}{q}))}$
			\\
			&{\footnotesize (Proposition \ref{prop:err:dyadic:Lq})}&{
				{\footnotesize\cite[Theorem 6.10]{BDSU14}, sparse grid}}\\
			\cmidrule{1-3}
			$L_{\infty}(\T^d)$  & $\lesssim M^{-\frac{\alpha-\frac{1}{2}}{2}}(\log
			M)^{\frac{d-2}{2}(\alpha-\frac{1}{2})+\frac{d-1}{2}}$ & $\asymp
			M^{-\alpha+\frac{1}{2}}(\log M)^{(d-1)\alpha}$ \\
			&{\footnotesize (Proposition \ref{thm:alias_err_infty:dyadic})}&{
				{\footnotesize\cite[Theorem 6.10]{BDSU14}, sparse grid}}\\
			\cmidrule{1-3}
			$\mathcal{H}^{\gamma}(\T^d)$&  $\lesssim M^{-\frac{\alpha-\gamma}{2}}(\log
			M)^{\frac{d-2}{2}(\alpha-\gamma)+\frac{d-1}{2}}$ & $\asymp
			M^{-(\alpha-\gamma)}$\\
			&{\footnotesize (Proposition \ref{prop:err:dyadic})}&{
				{\footnotesize\cite[Theorem 6.7]{BDSU14}, energy sparse grid}}\\
			\cmidrule{1-3}
			$\mathcal{H}^{\gamma}_{\mathrm{mix}}(\T^d)$&  $\lesssim
			M^{-\frac{\alpha-\gamma}{2}}(\log
			M)^{\frac{d-2}{2}(\alpha-\gamma)+\frac{d-1}{2}}$ & $\asymp
			M^{-(\alpha-\gamma)}(\log M)^{(d-1)(\alpha-\gamma)}$ \\
			&{\footnotesize (Theorem \ref{thm:alias_err:dyadic})}&{
				{\footnotesize\cite[Theorem 6.10]{BDSU14}, sparse grid}}\\
			\bottomrule
		\end{tabular}
	\end{center}
	\caption{Upper bounds of sampling numbers in the setting
		$\mathcal{H}^{\alpha}_{\text{mix}}(\T^d)\rightarrow Y$ for different sampling
		methods.
		Smoothness parameters are chosen from
		$\alpha>\max\{\gamma,\frac{1}{2}\}$, $\gamma>0$, and $2<q<\infty$. The upper bounds on
		$g_M^{\mathrm{latt}_1}$ are realized by the CBC rank-1
		lattice.}\label{tab:mainres:overview_general_r1l}
\end{table}
\begin{table}[h]
	\centering
	\begin{tabular}{lll}
		\toprule
		Y & $g_M^{\mathrm{latt}_1}(\mathcal{H}^{\alpha,\beta}(\T^d),Y)$ &
		$g_M^{\mathrm{lin}}(\mathcal{H}^{\alpha,\beta}(\T^d),Y)$ \\
		\cmidrule{1-3}
		$ L_2(\T^d)$  & $\lesssim M^{-\frac{\alpha+\beta}{2}}$ & $\asymp M^{-(\alpha+\beta)}$ \\
		&{\footnotesize \cite[Theorem 4.7]{KaPoVo14}}&{ {\footnotesize\cite[Theorem 6.10]{BDSU14}}}\\
		\cmidrule{1-3}
		$ L_q(\T^d)$  & $\lesssim M^{-\frac{\alpha-(\frac{1}{2}-\frac{1}{q})+\beta}{2}}(\log
		M)^{\frac{d-2}{2}(\alpha-(\frac{1}{2}-\frac{1}{q})+\beta)}$ &
		$\lesssim M^{-(\alpha-(\frac{1}{2}-\frac{1}{q})+\beta)}$ \\
		&{\footnotesize (Proposition \ref{prop:err:dyadic:Lq})}&{ {\footnotesize\cite[*]{BDSU14}}}\\
		\cmidrule{1-3}
		$L_{\infty}(\T^d)$  & $\lesssim M^{-\frac{\alpha+\beta-\frac{1}{2}}{2}}$ & $\lesssim M^{-(\alpha+\beta)+\frac{1}{2}}$ \\
		&{\footnotesize (Proposition \ref{thm:alias_err_infty:dyadic})}&
		{ {\footnotesize\cite[*]{BDSU14}}}\\
		\cmidrule{1-3}
		$\mathcal{H}^{\gamma}(\T^d)$&  $\lesssim M^{-\frac{\alpha+\beta-\gamma}{2}}(\log
		M)^{\frac{d-2}{2}(\alpha+\beta-\gamma)}$ & $\asymp
		M^{-(\alpha+\beta-\gamma)}$ \\
		&{\footnotesize (Proposition \ref{prop:err:dyadic})}&{ {\footnotesize\cite[Theorem 6.7]{BDSU14}}}\\
		\cmidrule{1-3}
		$\mathcal{H}^{\gamma}_{\mathrm{mix}}(\T^d)$&  $\lesssim M^{-\frac{\alpha+\beta-\gamma}{2}}(\log
		M)^{\frac{d-2}{2}(\alpha+\beta-\gamma)}$ & $\asymp
		M^{-(\alpha+\beta-\gamma)}$ \\
		&{\footnotesize (Theorem \ref{thm:alias_err:dyadic})}&{ {\footnotesize\cite[*]{BDSU14}}}\\
		\bottomrule
	\end{tabular}
	\caption{Upper bounds for sampling numbers for different sampling methods.
	Smoothness parameters are chosen from
	$\beta<0$, $\alpha+\beta>\max\{\gamma,\frac{1}{2}\}$, $\gamma>0$, and $2<q<\infty$. Best known bounds based on energy sparse grid sampling.
	References marked with $^*$ means that the result is not stated there explicitly but follows with the same method
	therein.}
	\label{tab:mainres:overview_general_r1l_2}
	\end{table}

	\begin{table}[h]
		\centering
		\begin{tabular}{lll}
			\toprule
			$Y$ & $g_M^{\mathrm{latt}_1}(\mathcal{H}^{\alpha}_{\text{mix}}(\T^2),Y)$ &
			$g_M^{\mathrm{lin}}(\mathcal{H}^{\alpha}_{\text{mix}}(\T^2),Y)$ \\
			\cmidrule{1-3}
			$L_2(\T^2)$  & $\asymp M^{-\frac{\alpha}{2}}$ & $\lesssim M^{-\alpha}(\log M)^{\alpha+\frac{1}{2}}$ \\
			&{\footnotesize (Theorem \ref{thm:alias_2d_err:dyadic})}&{ {\footnotesize\cite[Theorem 6.10]{BDSU14}, sparse grid}}\\
			\cmidrule{1-3}
			$L_{\infty}(\T^2)$  & $ \lesssim M^{-\frac{\alpha-\frac{1}{2}}{2}}$ & $\asymp M^{-\alpha+\frac{1}{2}}(\log M)^{\alpha}$
			\\
			&{\footnotesize (Proposition \ref{prop:err_2d_Linf})}&{ {\footnotesize\cite[Theorem 6.10]{BDSU14}, sparse grid}}\\
			\cmidrule{1-3}
			$\mathcal{H}^{\gamma}(\T^2)$&  $\asymp M^{-\frac{\alpha-\gamma}{2}}$ & $\asymp M^{-(\alpha-\gamma)}$ \\
			&{\footnotesize (Theorem \ref{thm:alias_2d_err:dyadic})}&{ {\footnotesize\cite[Theorem 6.7]{BDSU14}, energy sparse
					grid}}\\
			\cmidrule{1-3}
			$\mathcal{H}^{\gamma}_{\mathrm{mix}}(\T^2)$&  $\lesssim M^{-\frac{\alpha-\gamma}{2}}(\log M)^{\frac{1}{2}}$ & $\asymp
			M^{-(\alpha-\gamma)}(\log
			M)^{\alpha-\gamma}$ \\
			&{\footnotesize (Remark \ref{rem:err_2d_Hmix})}&{ {\footnotesize\cite[Theorem 6.10]{BDSU14}, sparse grid}}\\
			\bottomrule
		\end{tabular}
\caption{Upper bounds for sampling rates for different sampling methods.
Smoothness parameters are chosen from
$\alpha>\frac{1}{2}$, $\alpha>\gamma > 0$. The upper bounds for $g_M^{\mathrm{latt}_1}$ are realized either by the Fibonacci or CBC-generated lattice.
}
			\label{tab:mainres:overview_fibonacci_r1l}
		\end{table}
 		
{\bf Notation.} As usual, $\N$ denotes the natural numbers, $\N_0$ the non-negative integers,
$\mathbb{Z}$ the integers and $\mathbb{R}$ the real numbers. With $\T$ we denote the torus represented by the interval
$[0,1)$.
The letter $d$ is always reserved for the dimension in $\Z$, $\R$, $\N$, and $\T$.
For $0<p\leq \infty$ and $\zb x\in \R^d$ we denote $\|\zb x\|_p = (\sum_{i=1}^d |x_i|^p)^{1/p}$ with the
usual modification for $p=\infty$. The norm of an element $x \in X$ is denoted by $\|x|X\|$. 
If $X$ and $Y$ are two Banach spaces, the norm of an
operator $A\colon X \to Y$ will be denoted by $\|A|X\to Y\|$. The symbol $X \hookrightarrow Y$ indicates that there is a
continuous embedding from $X$ into $Y$. The relation $a_n\lesssim b_n$ means that there is a constant $c>0$ independent
of the context relevant
parameters such that $a_n\le c\,b_n$ for all $n$ belonging to a certain subset of $\N$, often $\N$ itself. We write $a_n
\asymp b_n$ if 
$a_n\lesssim b_n$ and $b_n \lesssim a_n$ holds. 

\section{Definitions and prerequisites}
The well known fact that decay properties of Fourier coefficients of a periodic function $f$ can be rephrased in smoothness properties of $f$ motivates to define the weighted Hilbert spaces
\begin{align} \label{def:hab}
\mathcal{H}^{\alpha,\beta}(\T^d):=\left\{f\in L_2(\T^d)\colon\|f|\mathcal{H}^{\alpha,\beta}(\T^d)\|^2:=\sum_{\zb k\in\Z^d}|\hat{f}_{\zb k}|^2(1+\|\zb k\|_2^2)^{\beta}\prod_{s=1}^d
(1+|k_s|^2)^{\alpha}<\infty\right\},
\end{align}
that mainly depend on the smoothness parameters $\alpha,\beta\in\R$, $\min\{\alpha,\alpha+\beta\}>0$. It is easy to show that for integer $\alpha,\beta\in \N_0$ these spaces coincide with the spaces defined in \eqref{eq:H_alpha_beta_integer}. Furthermore in case $\alpha=0$ and $\beta\geq 0$ these spaces coincide with isotropic Sobolev spaces, therefore we use the definition $\mathcal{H^{\beta}}(\T^d):=\mathcal{H}^{0,\beta}(\T^d)$. For $\alpha\geq0$ and $\beta=0$ the spaces $\mathcal{H}^{\alpha,0}(\T^d)$ coincide with $\mathcal{H}^{\alpha}_{\text{mix}}(\T^d)$, i.e. with the Sobolev spaces of dominating mixed smoothness, and we use the definition $\mathcal{H}^{\alpha}_{\text{mix}}(\T^d):=\mathcal{H}^{\alpha,0}(\T^d)$. 
Since we want to deal with sampling, we are interested in continuous functions. 
\begin{Lemma}\label{lem:sobolevembedding}
	Let $\alpha,\beta\in\R$ with $\min\{\alpha,\alpha+\beta\}>\frac{1}{2}$. Then
	$$\mathcal{H}^{\alpha,\beta}(\T^d)\hookrightarrow C(\T^d).$$
\end{Lemma}
\begin{proof}
	We refer to \cite[Theorem 2.9]{BDSU14}.
\end{proof}
The Fourier partial sum of a function $f\in L_1(\T^d)$ with respect to the frequency index set $I\subset\Z^d$, $|I|<\infty$, is defined by
$$S_If:=\sum_{\zb k\in I}\hat{f}_{\zb k}\e^{2\pi\ii\zb k\cdot\circ},$$
where 
\begin{align}
\hat{f}_{\zb k}:=\int_{\T^d}f(\zb x)\e^{-2\pi\ii\zb k\cdot\zb x}\mathrm{d}\zb x\label{eq:def_Fourier_coefficient}
\end{align}
are the usual Fourier coefficients of $f$.

We approximate the Fourier coefficients $\hat{f}_{\zb k}$, $\zb k\in I$, based on sampling values taken at the nodes of a rank-1 lattice 
\begin{align*}
\Lambda(\zb z, M):=\left\{\frac{j}{M}\zb z\bmod \zb 1\colon j=0,\ldots,M-1\right\}\subset\T^d,
\end{align*}
where $\zb z\in\Z^d$
is the generating vector and $M\in\N$ is the lattice size.
In particular, we apply the quasi-Monte Carlo rule defined by the rank-1 lattice $\Lambda(\zb z,M)$ on the integrand in \eqref{eq:def_Fourier_coefficient}, i.e.,
$$
\hat{f}_{\zb k}^{\Lambda(\zb z,M)}:=\frac{1}{M}\sum_{j=0}^{M-1}f\Big(\frac{j}{M}\zb z\Big)\e^{-2\pi\ii\frac{j}{M}\zb k\cdot\zb z}.
$$
Accordingly, we define the 
rank-1 lattice sampling operator $S_I^{\Lambda(\zb z,M)}$ by
\begin{align}
S_I^{\Lambda(\zb z,M)}f&:=\sum_{\zb k\in I}\hat{f}_{\zb k}^{\Lambda(\zb z,M)}\e^{2\pi\ii\zb k\cdot\circ}\label{eqn:sampling_operator}.
\end{align}

We call a rank-1 lattice $\Lambda(\zb z,M)$ \emph{reconstructing rank-1 lattice} for the frequency index set $I\subset\Z^d$, $|I|<\infty$, if the
sampling operator $S_I^{\Lambda(\zb z,M)}$ reproduces all trigonometric polynomials with frequencies supported on I, i.e.,
$S_I^{\Lambda(\zb z,M)}p=p$ holds for all trigonometric polynomials 
\begin{equation} p\in\Pi_I:=\sspan\{\e^{2\pi\ii\zb k\cdot\circ}\colon\zb k\in I\}\label{def:trigpol}.
\end{equation}
The condition
\begin{align}
\zb k^1\cdot\zb z&\not\equiv\zb k^2\cdot\zb z\imod{M}\quad\text{ for all }\zb k^1,\zb k^2\in I\text{, }\zb k^1\neq\zb k^2\text{,}\label{eqn:reco_prop1}
\end{align}
has to be fulfilled in order to guarantee that $\Lambda(\zb z,M)$ is a reconstructing rank-1 lattice for
the frequency index set $I$. One can show, that the condition in \eqref{eqn:reco_prop1} is not only sufficient but also necessary.
In the following sections, we frequently use the so-called difference set $\mathcal{D}(I)$ of a frequency index set $I\subset\Z^d$, $|I|<\infty$,
$$
\mathcal{D}(I):=\left\{\zb k\in\Z^d\colon\zb k=\zb h^1-\zb h^2,\,\zb h^1,\zb h^2\in  I\right\}.
$$
This definition allows for the reformulation of \eqref{eqn:reco_prop1} in terms of the difference set $\mathcal{D}(I)$, i.e.,
\begin{align}
 \zb k\cdot\zb z&\not\equiv0\imod{M}\quad\text{ for all }\zb k\in\mathcal{D}(I)\setminus\{\zb 0\}.\label{eqn:reco_prop2} 
\end{align}
Furthermore, we define the dual lattice
$$\Lambda(\zb z,M)^\perp:=\{\zb h\in\Z^d\colon\zb h\cdot\zb z\equiv 0\imod{M}\}$$
of the rank-1 lattice $\Lambda(\zb z,M)$. We use this definition in order to characterize
the reconstruction property of a rank-1 lattice $\Lambda (\zb z,M)$ for a frequency index set $I$.
A rank-1 lattice $\Lambda (\zb z,M)$ is a reconstructing rank-1 lattice for the frequency index set $I$,
$1\le|I|<\infty$, iff 
\begin{align}
\Lambda(\zb z,M)^\perp\cap\mathcal{D}(I)=\{\zb 0\}\label{eqn:dual_lattice_and_difference_set}
\end{align}
holds.
This means the conditions \eqref{eqn:reco_prop1}, \eqref{eqn:reco_prop2} and \eqref{eqn:dual_lattice_and_difference_set} are equivalent,
see also \cite{Kae2013}.
In order to approximate functions $f\in\mathcal{H}^{\alpha,\beta}(\T^d)$ using trigonometric polynomials,
we have to carefully choose the spaces $\Pi_I$ (cf. \eqref{def:trigpol}) of these trigonometric polynomials. Clearly, the spaces $\Pi_I$ are
described by the corresponding frequency index set $I$.
For technical reasons, we use so-called generalized dyadic hyperbolic crosses, 
\begin{align}
I=H_R^{d,T}:=\bigcup_{\zb j\in J_R^{d,T}}Q_{\zb j},\label{def:HRdT}
\end{align}
\begin{figure}
	\begin{minipage}{0.48\textwidth}
		\centering
		\begin{tikzpicture}[scale=0.09]
		\pgfmathsetmacro{\m}{5};
		\pgfmathsetmacro{\rf}{1};
		\pgfmathsetmacro{\rs}{1};
		\pgfmathtruncatemacro{\xl}{2^(\m / \rf)}
		\pgfmathtruncatemacro{\yl}{2^(\m / \rs)}
		\pgfmathtruncatemacro{\limf}{round(\m / \rf)}
		
		\pgfmathtruncatemacro{\lima}{round(\m / \rf)}
		\pgfmathtruncatemacro{\limb}{round(\m / \rs)}
		
		\foreach \x in {1,...,\lima}
		{
			\pgfmathtruncatemacro{\lim}{round((\m-\rf*\x)/(\rs))}
			\foreach \y in {1,...,\limb}
			{
				\pgfmathtruncatemacro{\xm}{\x-1}
				\pgfmathtruncatemacro{\ym}{\y-1}
				\draw[dashed] (2^\xm,2^\ym)-- (2^\x,2^\ym) --(2^\x,2^\y) -- (2^\xm,2^\y) --
				(2^\xm,2^\ym);
				\draw[dashed] (-2^\xm,-2^\ym)-- (-2^\x,-2^\ym) --(-2^\x,-2^\y) --
				(-2^\xm,-2^\y) -- (-2^\xm,-2^\ym) ;
				\draw[dashed] (-2^\xm,2^\ym)-- (-2^\x,2^\ym) --(-2^\x,2^\y) --
				(-2^\xm,2^\y) -- (-2^\xm,2^\ym) ;
				\draw[dashed] (2^\xm,-2^\ym)-- (2^\x,-2^\ym) --(2^\x,-2^\y) --
				(2^\xm,-2^\y) -- (2^\xm,-2^\ym) ;
			}
		}
		
		\foreach \x in {0,...,\limf}
		{
			\pgfmathtruncatemacro{\lim}{round((\m-\rf*\x)/(\rs))}
			\foreach \y in {0,...,\lim}
			{
				\pgfmathtruncatemacro{\xm}{\x-1}
				\pgfmathtruncatemacro{\ym}{\y-1}
				\ifnum \x<1 \relax%
				\ifnum \y<1 \relax%
				draw[color=black,very thick] (-1,-1)-- (1,-1) --(1,1) -- (-1,1) --
				(-1,-1);
				\else
				\draw[color=black,very thick] (-1,-2^\ym)-- (1,-2^\ym) --(1,-2^\y) --
				(-1,-2^\y) -- (-1,-2^\ym) ;
				\draw[color=black,very thick] (-1,2^\ym)-- (1,2^\ym) --(1,2^\y) --
				(-1,2^\y) -- (-1,2^\ym) ;
				\fi
				\else%
				\ifnum \y<1 \relax%
				\ifnum \x<1 \relax%
				
				\else
				\draw[color=black,very thick] (2^\xm,-1)-- (2^\x,-1) --(2^\x,1) --
				(2^\xm,1) -- (2^\xm,-1) ;
				\draw[color=black,very thick] (-2^\xm,-1)-- (-2^\x,-1) --(-2^\x,1)
				--
				(-2^\xm,1) -- (-2^\xm,-1);
				\fi
				\else
				\draw[color=black,very thick] (2^\xm,2^\ym)-- (2^\x,2^\ym) --(2^\x,2^\y) --
				(2^\xm,2^\y) -- (2^\xm,2^\ym);
				\draw[color=black,very thick] (-2^\xm,-2^\ym)-- (-2^\x,-2^\ym)
				--(-2^\x,-2^\y) -- (-2^\xm,-2^\y) -- (-2^\xm,-2^\ym) ;
				\draw[color=black,very thick] (-2^\xm,2^\ym)-- (-2^\x,2^\ym)
				--(-2^\x,2^\y) -- (-2^\xm,2^\y) -- (-2^\xm,2^\ym) ;
				\draw[color=black,very thick] (2^\xm,-2^\ym)-- (2^\x,-2^\ym)
				--(2^\x,-2^\y) -- (2^\xm,-2^\y) -- (2^\xm,-2^\ym) ;
				\fi
				\fi	
			}}
			\end{tikzpicture}
		\end{minipage}
		\begin{minipage}{0.48\textwidth}
			\centering
			\begin{tikzpicture}[scale=0.24]
			\begin{scope}[thick,font=\scriptsize]
			\draw [->] (0,0) -- (22,0) node [below] {$k_1$};
			\draw [->] (0,0) -- (0,22) node [above] {$k_2$};
			
			
			\end{scope}
			\node at (0,0) {$\bullet$};
			\node at (0,1) {$\bullet$};
			\node at (0,2) {$\bullet$};
			\node at (0,3) {$\bullet$};
			\node at (0,4) {$\bullet$};
			\node at (0,5) {$\bullet$};
			\node at (0,6) {$\bullet$};
			\node at (0,7) {$\bullet$};
			\node at (0,8) {$\bullet$};
			\node at (0,9) {$\bullet$};
			\node at (0,10) {$\bullet$};
			\node at (0,11) {$\bullet$};
			\node at (0,12) {$\bullet$};
			\node at (0,13) {$\bullet$};
			\node at (0,14) {$\bullet$};
			\node at (0,15) {$\bullet$};
			\node at (0,16) {$\bullet$};
			\node at (0,17) {$\bullet$};
			\node at (0,18) {$\bullet$};
			\node at (0,19) {$\bullet$};
			\node at (0,20) {$\bullet$};
			\node at (1,0) {$\bullet$};
			\node at (1,1) {$\bullet$};
			\node at (1,2) {$\bullet$};
			\node at (1,3) {$\bullet$};
			\node at (1,4) {$\bullet$};
			\node at (1,5) {$\bullet$};
			\node at (1,6) {$\bullet$};
			\node at (1,7) {$\bullet$};
			\node at (1,8) {$\bullet$};
			\node at (1,9) {$\bullet$};
			\node at (1,10) {$\bullet$};
			\node at (1,11) {$\bullet$};
			\node at (1,12) {$\bullet$};
			\node at (1,13) {$\bullet$};
			\node at (1,14) {$\bullet$};
			\node at (1,15) {$\bullet$};
			\node at (1,16) {$\bullet$};
			\node at (1,17) {$\bullet$};
			\node at (1,18) {$\bullet$};
			\node at (2,0) {$\bullet$};
			\node at (2,1) {$\bullet$};
			\node at (2,2) {$\bullet$};
			\node at (2,3) {$\bullet$};
			\node at (2,4) {$\bullet$};
			\node at (2,5) {$\bullet$};
			\node at (2,6) {$\bullet$};
			\node at (2,7) {$\bullet$};
			\node at (2,8) {$\bullet$};
			\node at (2,9) {$\bullet$};
			\node at (2,10) {$\bullet$};
			\node at (2,11) {$\bullet$};
			\node at (2,12) {$\bullet$};
			\node at (2,13) {$\bullet$};
			\node at (2,14) {$\bullet$};
			\node at (2,15) {$\bullet$};
			\node at (2,16) {$\bullet$};
			\node at (3,0) {$\bullet$};
			\node at (3,1) {$\bullet$};
			\node at (3,2) {$\bullet$};
			\node at (3,3) {$\bullet$};
			\node at (3,4) {$\bullet$};
			\node at (3,5) {$\bullet$};
			\node at (3,6) {$\bullet$};
			\node at (3,7) {$\bullet$};
			\node at (3,8) {$\bullet$};
			\node at (3,9) {$\bullet$};
			\node at (3,10) {$\bullet$};
			\node at (3,11) {$\bullet$};
			\node at (3,12) {$\bullet$};
			\node at (3,13) {$\bullet$};
			\node at (3,14) {$\bullet$};
			\node at (4,0) {$\bullet$};
			\node at (4,1) {$\bullet$};
			\node at (4,2) {$\bullet$};
			\node at (4,3) {$\bullet$};
			\node at (4,4) {$\bullet$};
			\node at (4,5) {$\bullet$};
			\node at (4,6) {$\bullet$};
			\node at (4,7) {$\bullet$};
			\node at (4,8) {$\bullet$};
			\node at (4,9) {$\bullet$};
			\node at (4,10) {$\bullet$};
			\node at (4,11) {$\bullet$};
			\node at (4,12) {$\bullet$};
			\node at (5,0) {$\bullet$};
			\node at (5,1) {$\bullet$};
			\node at (5,2) {$\bullet$};
			\node at (5,3) {$\bullet$};
			\node at (5,4) {$\bullet$};
			\node at (5,5) {$\bullet$};
			\node at (5,6) {$\bullet$};
			\node at (5,7) {$\bullet$};
			\node at (5,8) {$\bullet$};
			\node at (5,9) {$\bullet$};
			\node at (5,10) {$\bullet$};
			\node at (6,0) {$\bullet$};
			\node at (6,1) {$\bullet$};
			\node at (6,2) {$\bullet$};
			\node at (6,3) {$\bullet$};
			\node at (6,4) {$\bullet$};
			\node at (6,5) {$\bullet$};
			\node at (6,6) {$\bullet$};
			\node at (6,7) {$\bullet$};
			\node at (6,8) {$\bullet$};
			\node at (7,0) {$\bullet$};
			\node at (7,1) {$\bullet$};
			\node at (7,2) {$\bullet$};
			\node at (7,3) {$\bullet$};
			\node at (7,4) {$\bullet$};
			\node at (7,4) {$\bullet$};
			\node at (7,5) {$\bullet$};
			\node at (7,6) {$\bullet$};
			\node at (8,0) {$\bullet$};
			\node at (8,1) {$\bullet$};
			\node at (8,2) {$\bullet$};
			\node at (8,3) {$\bullet$};
			\node at (8,4) {$\bullet$};
			\node at (8,5) {$\bullet$};
			\node at (8,6) {$\bullet$};
			\node at (9,0) {$\bullet$};
			\node at (9,1) {$\bullet$};
			\node at (9,2) {$\bullet$};
			\node at (9,3) {$\bullet$};
			\node at (9,4) {$\bullet$};
			\node at (9,5) {$\bullet$};
			\node at (10,0) {$\bullet$};
			\node at (10,1) {$\bullet$};
			\node at (10,2) {$\bullet$};
			\node at (10,3) {$\bullet$};
			\node at (10,4) {$\bullet$};
			\node at (10,5) {$\bullet$};
			\node at (11,0) {$\bullet$};
			\node at (11,1) {$\bullet$};
			\node at (11,2) {$\bullet$};
			\node at (11,3) {$\bullet$};
			\node at (11,4) {$\bullet$};
			\node at (12,0) {$\bullet$};
			\node at (12,1) {$\bullet$};
			\node at (12,2) {$\bullet$};
			\node at (12,3) {$\bullet$};
			\node at (12,4) {$\bullet$};
			\node at (13,0) {$\bullet$};
			\node at (13,1) {$\bullet$};
			\node at (13,2) {$\bullet$};
			\node at (13,3) {$\bullet$};
			\node at (14,0) {$\bullet$};
			\node at (14,1) {$\bullet$};
			\node at (14,2) {$\bullet$};
			\node at (14,3) {$\bullet$};
			\node at (15,0) {$\bullet$};
			\node at (15,1) {$\bullet$};
			\node at (15,2) {$\bullet$};
			\node at (16,0) {$\bullet$};
			\node at (16,1) {$\bullet$};
			\node at (16,2) {$\bullet$};
			\node at (17,0) {$\bullet$};
			\node at (17,1) {$\bullet$};
			\node at (18,0) {$\bullet$};
			\node at (18,1) {$\bullet$};
			\node at (19,0) {$\bullet$};
			\node at (20,0) {$\bullet$};
			
			\end{tikzpicture}
		\end{minipage}
		\caption{$H^{2,0}_4$ and $J^{2,0.5}_{20}$}\label{fig:hypcross}
	\end{figure}
cf. Figure \ref{fig:hypcross}, where $R\ge1$ denotes the refinement, $T\in[0,1)$ is an additional parameter,
$$J_R^{d,T}:=\{\zb j\in\N_0^d\colon\|\zb j\|_1-T\|\zb j\|_\infty\le (1-T)R+d-1\},$$
and $Q_{\zb j}:=\bigtimes_{s=1}^d Q_{j_s}$ are sets of tensorized dyadic intervals
\begin{align}
Q_{j}&:=
\begin{cases}
\{-1,0,1\}&:\;j=0,\\
([-2^{j},-2^{j-1}-1]\cup[2^{j-1}+1,2^{j}])\cap\Z&:\;j>0,
\end{cases}\label{equ:def_Q_j}
\end{align}
cf. \cite{KnDiss}. 
\begin{Lemma}\label{lem:card_HRdT}
	Let the dimension $d\in\N$, the parameter $T\in[0,1)$, and the refinement $R\ge 1$, be given.
	Then, we estimate the cardinality of the index set  $H_R^{d,T}$ by
	$$
	|H_R^{d,T}|\asymp
	\begin{cases}
	2^RR^{d-1}&:\;T=0,\\
	2^R&:\;0<T<1.
	\end{cases}
	$$
\end{Lemma}
\begin{proof}
	The assertion for the upper bound follows directly from \cite[Lemma 4.2]{GrKn09}. For a proof including the lower bound we refer to \cite[Lemma 6.6]{BDSU14}. 
\end{proof}

Having fixed the index set $I=H^{d,T}_R$ an important question is the existence of a reconstructing lattice for it. If there is such a lattice, out of how many points does it consist? Can we explicitly
construct it? The following lemma answers these questions.
	\begin{Lemma}\label{lem:universal_bounds_M}
		Let the parameters $T\in[0,1)$, $R\ge1$, and the dimension $d\in\N$, $d\ge2$, be given.
		Then, there exists
		a reconstructing rank-1 lattice $\Lambda(\zb z,M)$ for $H_R^{d,T}$ which fulfills
		\begin{align*}
		2^{2R-2}\le M\lesssim
		\begin{cases}
		2^{2R}&:\;T>0,\\
		2^{2R}R^{d-2}&:\;T=0. 
		\end{cases}
		\end{align*}
		Moreover, each reconstructing rank-1 lattice $\Lambda(\zb z,M)$ for $H_R^{d,T}$ fulfills the lower bound.
	\end{Lemma}
	\begin{proof}
		For $T=0$, a detailed proof of the bounds can be found in \cite{Kae2012}.
		In the case $T\in(0,1)$, one proves the lower bound using the same way as used for $T=0$.
		The corresponding upper bound follows directly from \cite[Cor. 1]{Kae2013} and 
		$H_R^{d,T}\subset [-|H_R^{d,T}|,|H_R^{d,T}|]^d$ and $|H_R^{d,T}|\lesssim 2^R$.
	\end{proof}
	A lattice fulfilling these properties can be explicitly constructed using a component-by-component (CBC) optimization strategy for the generating vector $\zb z$. For more details on that
algorithm we refer to \cite[Ch. 3]{kaemmererdiss}.

\section{Lower bounds and non-optimality}\label{sec:lower_bounds}
\begin{figure}[tb]
\centering
	\subfloat[\label{fig:ac}$\AC_{8}^{3}$]{\includegraphics[scale=0.98]{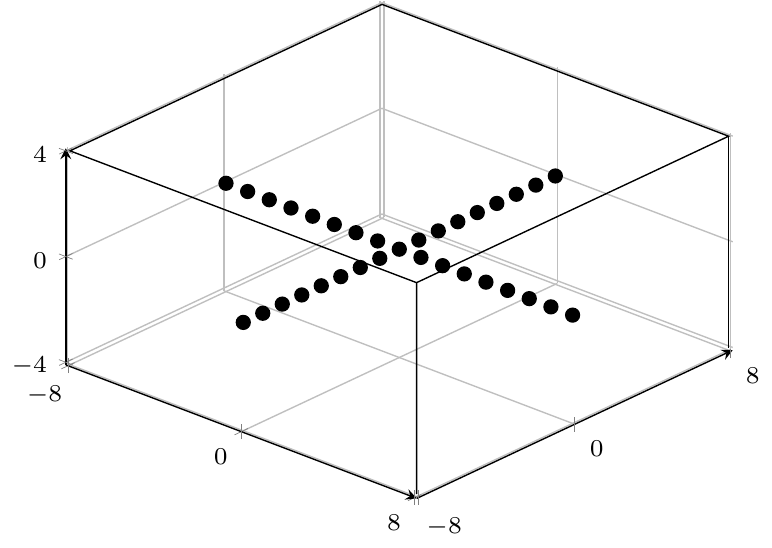}}
	\subfloat[\label{fig:acd}$P_{8}^2\setminus\{\zb 0\}$]{\includegraphics[scale=0.98]{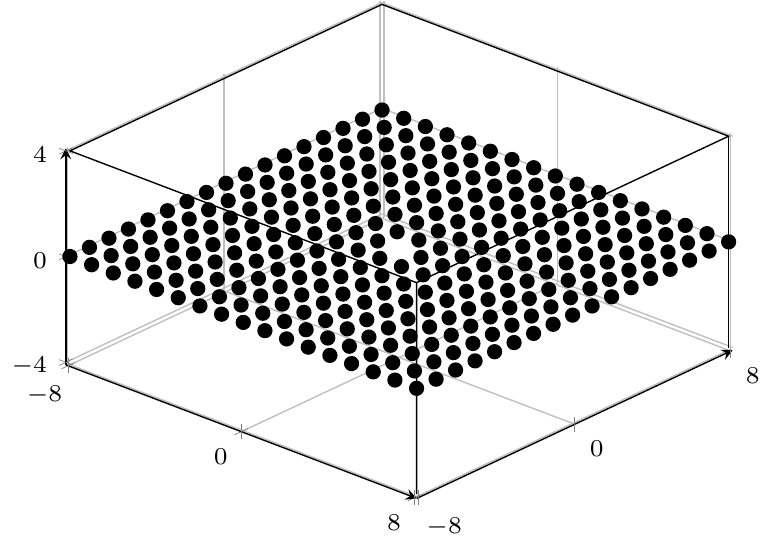}}
	\caption{Axis cross and subset of the difference set of the corresponding axis cross.}
\end{figure}

\begin{sloppypar}
	In this chapter we study lower bounds for the rank-1 lattice sampling numbers $g^{\mathrm{latt}_1}_M(\mathcal{H}^{\alpha,\beta}(\T^d),\mathcal{H}^{\gamma}(\T^d))$ and 
$g^{\mathrm{latt}_1}_M(\mathcal{H}^{\alpha,\beta}(\T^d),\mathcal{H}^{\gamma}_{\mathrm{mix}}(\T^d))$.
	At first we show, that each rank-1 lattice $\Lambda(\zb z,M)$, $\zb z\in\Z^d$, $d\ge 2$, and $M\in\N$,
	has at least one aliasing pair of frequency indices $\zb k^1,\zb k^2$ within the two-dimensional axis cross 
	$$
	\AC_{\sqrt{M}}^{d}:=\{\zb h\in\Z^2\times\underbrace{\{0\}\times\ldots\times\{0\}}_{d-2 \text{ times}}\colon\|\zb h\|_1=\|\zb h\|_\infty\le\sqrt{M}\}.
	$$
	For illustration, we depict $\AC_{8}^{3}$ in Figure \ref{fig:ac}.
	We can even show a more general result.
\end{sloppypar}

\begin{Lemma}\label{lem:aliasing_axis_cross}
	Let $\mathcal{X}:=\{\zb x_j\in\T^d\colon j=0,\ldots,M-1\}$, $d\ge 2$, be a sampling set of cardinality $|\mathcal{X}|=M$.
	In addition, we assume that
	\begin{align*}
	\sum_{j=0}^{M-1}\e^{2\pi\ii\zb k\cdot\zb x_j}\in\{0,M\}
	\text{ for all }
	\zb k\in P_{\sqrt{M}}^d:=\{-\floor{\sqrt{M}},\ldots,\floor{\sqrt{M}}\}^2\times\underbrace{\{0\}\times\ldots\times\{0\}}_{d-2 \text{ times}}.
	\end{align*}
	Then there exist at least two distinct indices $\zb k^1,\zb k^2\in\AC_{\sqrt{M}}^{d}$ within the axis cross $\AC_{\sqrt{M}}^{d}$ such that
	$\e^{2\pi\ii\zb k^1\cdot\zb x_j}=\e^{2\pi\ii\zb k^2\cdot\zb x_j}$  for all $j=0,\ldots,M-1$.
\end{Lemma}

\begin{proof}
	First, we assume 
	\begin{align}
	\sum_{j=0}^{M-1}\e^{2\pi\ii\zb h\cdot\zb x_j}=0 \text{ for all } \zb h\in P_{\sqrt{M}}^d\setminus\{\zb 0\},\label{eqn:assumption_proof}
	\end{align}
	cf. Figure \ref{fig:acd} for an illustration of the index set. 
	Consequently, for all $\zb h^1,\zb h^2\in \tilde{P}_{\sqrt{M}}^d:=\{0,\ldots,\floor{\sqrt{M}}\}^2\times\underbrace{\{0\}\times\ldots\times\{0\}}_{d-2 \text{ times}}$
	we achieve $\zb h^2-\zb h^1\in P_{\sqrt{M}}^d$ and
	\begin{align*}
	\sum_{j=0}^{M-1}\e^{2\pi\ii(\zb h^2-\zb h^1)\cdot\zb x_j}&=\begin{cases}M&:\;\zb h^2-\zb h^1=0\\0&\text{otherwise.}\end{cases}
	\end{align*}
	In matrix vector notation this means
	$$
	\zb A^*\zb A=M\zb I,
	$$
	where the matrix $\zb A=\Big(\e^{2\pi\ii\zb h\cdot\zb x_j}\Big)_{j=0,\ldots,M-1,\,\zb h\in \tilde{P}_{\sqrt{M}}^d}\in\C^{M\times\left(\floor{\sqrt{M}}+1\right)^2}$ must have full column rank.
	However, this is not possible due to the inequality $M<\left(\floor{\sqrt{M}}+1\right)^2$. Thus, the assumption given in \eqref{eqn:assumption_proof} does not hold in any case.
	\newline
	Accordingly, we consider the case where $\sum_{j=0}^{M-1}\e^{2\pi\ii\zb{h'}\cdot\zb x_j}=M$ for at least one 
	$\zb{h'}\in P_{\sqrt{M}}^d\setminus\{\zb 0\}$.
	Consequently, we observe $\e^{2\pi\ii\zb{h'}\cdot\zb x_j}=1$ for all $j=0,\ldots,M-1$.
	Then, for the frequency indices $\zb k^1=(h_1',0\ldots,0)^\top\in\AC_{\sqrt{M}}^{d}$ and $\zb k^2=(0,-h_2',0\ldots,0)^\top\in\AC_{\sqrt{M}}^{d}$, 
	the equalities $\e^{2\pi\ii\zb k^1\cdot\zb x_j}=\e^{2\pi\ii\zb k^2\cdot\zb x_j}$, $j=0,\ldots,M-1$, hold.
\end{proof}

As a consequence of the last considerations, we know that for each $d$-dimensional rank-1 lattice of size $M$, $d\ge2$,
there is at least one pair $\zb k^1,\zb k^2\in \AC_{\floor{\sqrt{M}}}^{d}=\AC_{\sqrt{M}}^{d}$ of frequencies within the two-dimensional axis cross of size $\sqrt{M}$ fulfilling
$$\zb k^1\cdot \zb z \equiv \zb k^2 \cdot z \imod M.$$
We call such a pair aliasing pair. As a consequence, we estimate the error of rank-1 lattice sampling operators from below as follows.

\begin{Theorem}\label{thm:lower_bound_omega_alpha_beta}
	Let the smoothness parameters $\alpha,\beta,\gamma\in\R$, $\alpha>\gamma-\beta \ge0$, $\alpha+\beta>\frac{1}{2}$.
	Then, we obtain
	\begin{align}
	g^{\mathrm{latt}_1}_M(\mathcal{H}^{\alpha,\beta}(\T^d),\mathcal{H}^{\gamma}(\T^d))
	&\ge 2^{-(\alpha+\beta-\gamma+1)/2} M^{-(\alpha+\beta-\gamma)/2} \label{eq:lowerbound_H_isotropic}
	\intertext{and}
	g^{\mathrm{latt}_1}_M(\mathcal{H}^{\alpha,\beta}(\T^d), \mathcal{H}^{\gamma}_{\mathrm{mix}}(\T^d))
	&\ge 2^{-(\alpha+\beta-\gamma+1)/2} M^{-(\alpha+\beta-\gamma)/2}.\label{eq:lowerbound_H_mix}
	\end{align}
	for all $M\in \N$.
\end{Theorem}
\begin{proof}
	For a given rank-1 lattice $\Lambda(\zb z,M)$, we construct the fooling function $\tilde{g}(\zb x):=\e^{2\pi\ii\zb k^1\cdot \zb x}-\e^{2\pi\ii\zb k^2\cdot \zb x}$,
	where  $\zb k^1, \zb k^2\in\AC_{\sqrt{M}}^{d}$ are aliasing frequency indices with respect to $\Lambda(\zb z,M)$, i.e., $\zb k^1\cdot\zb z\equiv\zb k^2\cdot\zb z\imod{M}$.
	These aliasing frequency indices exist due to Lemma \ref{lem:aliasing_axis_cross}. Using the notation
	$$\omega^{d,\alpha,\beta}(\zb k)^2:=\Big[\prod_{s=1}^d(1+|k_s|^2)\Big]^{\alpha}(1+\|\zb k\|_2^2)^{\beta},$$
	 the normalization of $\tilde{g}$ in $\mathcal{H}^{\alpha,\beta}(\T)$ is given by
	$$g(x):=\frac{\e^{2\pi\ii\zb k^1\cdot \zb x}-\e^{2\pi\ii\zb k^2\cdot \zb x}}{\sqrt{\omega^{d,\alpha,\beta}(\zb k^1)^2+\omega^{d,\alpha,\beta}(\zb k^2)^2}}.$$
	According to Lemma \ref{lem:aliasing_axis_cross}, the fooling function
	$g$ is zero at all sampling nodes  $\zb x_j\in\Lambda(\zb z,M)$ and
	we obtain
	\begin{align}
	\|g|\mathcal{H}^{\gamma}(\T^d)\|
	&=\frac{\sqrt{\omega^{d,0,\gamma}(\zb k^1)^2+\omega^{d,0,\gamma}(\zb k^2)^2}}{\sqrt{\omega^{d,\alpha,\beta}(\zb k^1)^2+\omega^{d,\alpha,\beta}(\zb k^2)^2}}.\nonumber 
	\intertext{W.l.o.g.  we assume $\|\zb k^1\|_\infty\ge\|\zb k^2\|_\infty$, i.e., $\omega^{d,0,\gamma}(\zb k^1)\ge\omega^{d,0,\gamma}(\zb k^2)$ and 
	$\omega^{d,\alpha,\beta}(\zb k^1)\ge\omega^{d,\alpha,\beta}(\zb k^2)$. We achieve}
		\|g|\mathcal{H}^{\gamma}(\T^d)\|
		&\geq \frac{\sqrt{\omega^{d,0,\gamma}(\zb k^1)^2}}{\sqrt{2\omega^{d,\alpha,\beta}(\zb k^1)^2}}=\frac{1}{\sqrt{2}\omega^{d,\alpha,\beta-\gamma}(\zb k^1)}.\label{eq:insertlower}
	\end{align}
	For $\zb k\in \AC_{\sqrt{M}}^{d}$ with $|k_1|=\|\zb k\|_{\infty}$ and $M\in\N$ we have
	\begin{align*}
		\omega^{d,\alpha,\beta-\gamma}(\zb k) &= (1+|k_1|^2)^{(\alpha+\beta-\gamma)/2}\\
		&\leq (1+M)^{(\alpha+\beta-\gamma)/2} \leq (2M)^{(\alpha+\beta-\gamma)/2}.
	\end{align*}
Inserting this into \eqref{eq:insertlower} yields 
\begin{align*}
		\|g|\mathcal{H}^{\gamma}(\T^d)\|
		&\geq 2^{-(\alpha+\beta-\gamma+1)/2}M^{-(\alpha+\beta-\gamma)/2}
\end{align*}
Now \eqref{eq:lowerbound_H_isotropic} follows by a standard argument. Let $A:\C^M\mapsto \mathcal{H}^{\gamma}(\T^d)$ be an arbitrary algorithm applied to $\Bigg(f(\zb 0),f\Big(\frac{1}{M}\zb
z\Big),\ldots,f\Big(\frac{M-1}{M}\zb z\Big)\Bigg)=\zb 0$. We estimate as follows
\begin{align*}
	 2^{-(\alpha+\beta-\gamma+1)/2}M^{-(\alpha+\beta-\gamma)/2}\leq& \|g|\mathcal{H}^{\gamma}(\T)\|\leq \frac{1}{2}(\|g-A(\zb 0)|\mathcal{H}^{\gamma}(\T)\|+\|-g-A(\zb
0)|{\mathcal{H}^{\gamma}(\T)})\|\\
	\leq& \max\{\|g-A(\zb 0)|{\mathcal{H}^{\gamma}(\T)}\|,\|-g-A(\zb 0)|{\mathcal{H}^{\gamma}(\T)}\|\}.
\end{align*}
Accordingly, each algorithm $A$ badly approximates at least one of the functions $g$ or $-g$.
Thus, we observe an infimum over the worst case errors of all algorithms $A$ 
\begin{align*}
	\operatorname{Samp}_{\Lambda(\zb z,M)}(\mathcal{H}^{\alpha,\beta}(\T^d),H^{\gamma}(\T^d)) \geq 2^{-(\alpha+\beta-\gamma+1)/2}M^{-(\alpha+\beta-\gamma)/2}.
\end{align*}
Finally the infimum over all rank-1 lattices with $M$ points yields
\begin{align*}
\operatorname{g}_{M}^{\mathrm{latt}_1}(\mathcal{H}^{\alpha,\beta}(\T^d),H^{\gamma}(\T^d)) \geq 2^{-(\alpha+\beta-\gamma+1)/2}M^{-(\alpha+\beta-\gamma)/2}.
\end{align*}
	The assertion in \eqref{eq:lowerbound_H_mix} can be proven analogously.
\end{proof}

Following attentively the last proof we recognize that the condition $\alpha+\beta>\frac{1}{2}$ plays no fundamental role in the estimations there. It is required for a well interpretation of the function evaluations in the definition of $g^{\mathrm{latt}_1}_M(\mathcal{H}^{\alpha,\beta}(\T^d),Y)$, which is given for continuous functions (cf. Lemma \ref{lem:sobolevembedding}). For $\min\{\alpha,\alpha+\beta\}>0$, a generalization of the last theorem can be achieved using the space
$$\mathcal{H}^{\alpha,\beta}(\T^d)\cap^* C(\T^d):=\left\{f\in C(\T^d)\colon\|f|\mathcal{H}^{\alpha,\beta}(\T^d)\|<\infty\right\},$$
equipped with the norm of $\mathcal{H}^{\alpha,\beta}(\T^d)$,
see \eqref{def:hab} for comparison.
Then the proof of Theorem \ref{thm:lower_bound_omega_alpha_beta} yields the following proposition.
\begin{Proposition}\label{prop:lower_bound_r1l_general}
	Let the smoothness parameters $\alpha,\beta,\gamma\in\R$, $\alpha>\gamma-\beta \ge0$, $\alpha+\beta>0$.
	Then, we obtain
	\begin{align*}
	g^{\mathrm{latt}_1}_M(\mathcal{H}^{\alpha,\beta}(\T^d)\cap^* C(\T^d),\mathcal{H}^{\gamma}(\T^d))
	&\ge 2^{-(\alpha+\beta-\gamma+1)/2} M^{-(\alpha+\beta-\gamma)/2}
	\intertext{and}
	g^{\mathrm{latt}_1}_M(\mathcal{H}^{\alpha,\beta}(\T^d)\cap^* C(\T^d), \mathcal{H}^{\gamma}_{\mathrm{mix}}(\T^d))
	&\ge 2^{-(\alpha+\beta-\gamma+1)/2} M^{-(\alpha+\beta-\gamma)/2}.
	\end{align*}
	for all $M\in \N$.
\end{Proposition}
\begin{Remark}
	We stress on the fact that even each $d$-dimensional rank-$s$ lattice of size $M$, where $d\ge2$ and $s\in\N$, $s\le d$,
	fulfills the requirements of Lemma \ref{lem:aliasing_axis_cross}, cf. \cite[Lemma 2.7]{SlJo94}.
	Consequently, there exists at least one aliasing pair
	$\zb k^1,\zb k^2\in\AC_{\sqrt{M}}^{d}$
	within the two-dimensional axis cross of size $\sqrt{M}$.
	This means we obtain the statements of Theorem \ref{thm:lower_bound_omega_alpha_beta}
	using the identical proof strategy.
\end{Remark}

%

\section{Improved upper bounds for $d>2$}\label{sec:upper_bounds}
In this section we study upper bounds for $g_M^{\mathrm{latt}_1}$. To do this, we consider approximation error estimates for
$S_{H^{d,T}_R}^{\Lambda(\zb z,M)}f$. To obtain these estimates the cardinality of the dual lattice $\Lambda(\zb z,M)^{\perp}$ intersected with rectangular boxes $\Omega$ plays an important role.
\begin{Lemma}\label{lem:chi_estimate}
Let $\Lambda(\zb z,M)$ be a rank-1 lattice generated by $\zb z \in \Z^d$ with $M$ points. Assume that the dual lattice $\Lambda(\zb z,M)^{\perp}$ is located outside the hyperbolic cross $H^{d,0}_R$, $R\geq1$, i.e.,
\begin{equation}\Lambda(\zb z,M)^{\perp}\cap H_R^{d,0}=\{\zb 0\}.\label{eq:341}\end{equation}
Then we have
\begin{equation}
|\Lambda(\zb z,M)^{\perp}\cap \Omega|\leq \begin{cases}2^{d+1}\frac{\vol \Omega}{2^R}&:\;\vol \Omega > 2^{R-1},\\
1&:\:\vol \Omega \leq 2^{R-1},
\end{cases}\label{eq:342}
\end{equation}
where $\Omega$ is an arbitrary rectangle with side-lenghts $\geq 1$.
\end{Lemma}
\begin{proof}
For two arbitrary distinct dual lattice points $\boldk^1,\boldk^2\in \Lambda(\zb z,M)^{\perp}$, $\boldk^1\neq \boldk^2$, we obtain $\boldk:=\boldk^1-\boldk^2\in \Lambda(\zb z,M)^{\perp}\setminus\{\zb 0\}$. As a consequence of \eqref{def:HRdT} and \eqref{eq:341}, the vector $\boldk$ belongs to a cuboid $Q_{\zb j}$ with $\|\zb j\|_1>R+d-1$. We achieve
$$\prod_{s=1}^d \max\{|k_s|,1\}
=\prod_{\substack{s=1\\j_s> 0}}^d|k_s|\geq
\prod_{\substack{s=1\\j_s> 0}}^d(2^{j_s-1}+1)
>
\prod_{\substack{s=1\\j_s> 0}}^d2^{j_s-1}
\geq2^{\|\zb j\|_1-d}
> 2^{R-1}.$$

{\em{Step 1.}} We prove the second case in \eqref{eq:342} by contradiction.
For any rectangle $\Omega:=[a_1,a_1+b_1]\times \ldots \times [a_d,a_d+b_d]$ with side lengths $b_s\geq 1$, $s=1,\ldots,d$, and  $\vol \Omega=\prod_{s=1}^db_s\leq 2^{R-1}$ we assume $|\Lambda(\zb z,M)^{\perp}\cap \Omega|\geq 2$ and we choose $\boldk^1,\boldk^2\in\Omega\cap \Lambda(\zb z,M)^{\perp}$, $\boldk^1\neq \boldk^2$.
Consequently, there is a $d$-dimensional cuboid $K\subset\Omega$ of side lengths $\geq1$ which contains the minimal cuboid with edges $\boldk^1$ and $\boldk^2$.
The volume of $K$ is bounded from below by
$\prod_{s=1}^d \max\{|k_s|,1\}> 2^{R-1}$, and hence larger than the volume of $\Omega$, which is in contradiction to the 
relation $K\subset\Omega$. Accordingly, there can not be more than one element within $\Lambda(\zb z,M)^{\perp}\cap \Omega$.

{\em{Step 2.}} We prove the first case and assume that $\Omega$ has volume larger than $2^{R-1}$. The sidelengths of $\Omega$ are denoted by $b_s$, $s=1,\ldots,d$. We construct
a disjoint covering/packing of $\Omega$ consisting of half side opened cuboids $B$ with sidelength $\ell_1,\ldots,\ell_d$ such
that $\ell_s\leq \max(1,b_s)$, $s=1,\ldots,d$, and $\operatorname{vol} B = 2^{R-1}$, cf. Figure \ref{fig:illustrate_Lemma} for illustration.
We need at most $2^d\frac{\vol \Omega}{2^{R-1}}$ of the cuboids $B$ in order to cover the set $\Omega$.
Due to Step 1, each $B$ contains at most one element from $\Lambda(\zb z,M)^{\perp}$.
Accordingly, the number of elements in $\Lambda(\zb z,M)^{\perp}\cap \Omega$ is bounded from above by
$2^{d+1}\frac{\vol \Omega}{2^R}$.
\begin{figure}[tb]
\centering
\begin{tikzpicture}[scale=1]
\foreach \x in {0,1.5,3,4.5,6,7.5} 
\foreach \y in {0,1,2,3} 
\draw[very thick,color=lightgray] (\x,\y) -- (\x+1.5,\y) -- (\x+1.5,1+\y) -- (\x,1+\y) --  (\x,0+\y);
\draw[very thick] (0,0) -- (8,0) -- (8,3.5) -- (0,3.5) -- (0,0);
\node at (0.6,1.2)[above] {$\ell_1$};
\node at (1.7,0.6)[right] {$\ell_2$};
\node at (0.75,0.5) {$B$};
\node at (4,-0.5) [below] {$b_1$};
\node at (-0.5,1.75)[left] {$b_2$};
\node at (4,2) {$\Omega$};
\draw[<->] (0,1.2) -- (1.5,1.2);
\draw[<->] (1.7,0) -- (1.7,1);
\draw[<->] (-0.5,0) -- (-0.5,3.5);
\draw[<->] (0,-0.5) -- (8,-0.5);

\end{tikzpicture}
\caption{The counting argument in Lemma \ref{lem:chi_estimate}.}\label{fig:illustrate_Lemma}
\end{figure}
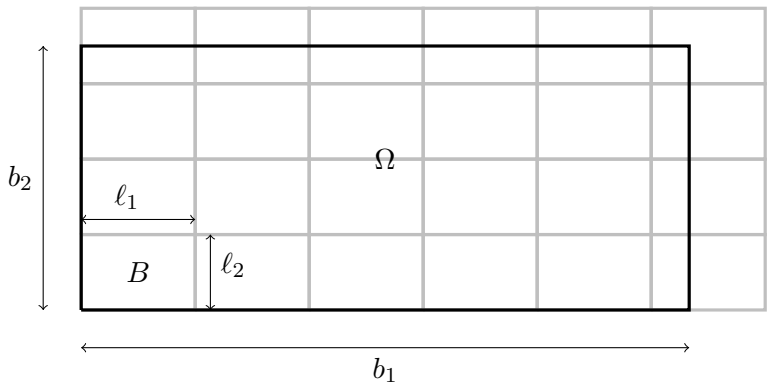
\end{proof}
\begin{Lemma}\label{lem:theta}
	Let the smoothness parameters $\alpha,\beta\in\R$, $\beta\le0$, $\alpha+\beta>1/2$, the refinement $R\ge 1$, and the parameter $T:=-\beta/\alpha$ be given. In addition, we assume that the rank-1 lattice $\Lambda(\zb z,M)$ is a reconstructing rank-1 lattice for the hyperbolic cross $H_R^{d,0}$. We define
	\begin{align}
	\theta^2_{\alpha,\beta}(\zb k,\zb z,M):=
	\sum_{\substack{\boldh\in \Lambda(\zb z,M)^{\perp}\\\boldh\neq \boldzero}}(1+\|\boldk+\boldh\|_2^2)^{-\beta}\prod_{s=1}^d(1+|k_s+h_s|^2)^{-\alpha}.
	\label{def:theta}
	\end{align}
	Then the estimate
	$$
	\theta^2_{\alpha,\beta}(\zb k,\zb z,M)\lesssim
	\begin{cases}
	2^{-2(\alpha+\beta)R} & :\;T>0,\\
	2^{-2\alpha R}R^{d-1} & :\;T=\beta=0
	\end{cases}
	$$
	holds for all $\zb k\in H_R^{d,0}$.
\end{Lemma}
\begin{proof}
	For $\zb k\in\Z^d$ and $\zb j\in \N_0^d$ we define the indicator function
	\begin{align*}
	\varphi_{\zb j}(\zb k)&:=\begin{cases}0 &:\;\zb k\not\in Q_{\zb j},\\1 &:\;\zb k\in Q_{\zb j},\end{cases}
	\end{align*}
	where $Q_{\zb j}$ is defined in \eqref{equ:def_Q_j}.
	We fix  $\zb k\in H_R^{d,0}$ and decompose the sum in \eqref{def:theta}, which yields
	\begin{align*}
	\theta^2_{\alpha,\beta}(\zb k,\zb z,M)&=
	\sum_{\substack{\boldh\in \Lambda(\zb z,M)^{\perp}\\\boldh\neq \boldzero}}\sum_{\zb j\in\N_0^d}\varphi_{\zb j}(\zb k+\zb h)(1+\|\boldk+\boldh\|_2^2)^{-\beta}\prod_{s=1}^d(1+|k_s+h_s|^2)^{-\alpha}.
	\end{align*}
	Since $\Lambda(\zb z,M)$ is a reconstructing rank-1 lattice for $H_R^{d,0}$, we know from \eqref{eqn:dual_lattice_and_difference_set} that 
	$$
	  \mathcal{D}(H_R^{d,0})\cap\left(\Lambda(\zb z,M)^\perp\setminus\{\zb 0\}\right)=\emptyset\,.
	$$
	This yields  $$\zb k^1+\zb h^1\neq \zb k^2+ \zb h^2$$ for all $\zb k^1,\zb k^2\in H^{d,0}_R$, $\zb k^1 \neq \zb k^2$, and $\zb h^1,\zb h^2\in \Lambda(\zb z,M)^{\perp}$ since otherwise 
	$\zb 0\neq \zb k^1- \zb k^2 = \zb h^2-\zb h^1 \in \Lambda(\zb z,M)^{\perp}$ which is in contradiction to \eqref{eqn:dual_lattice_and_difference_set}.
	In particular, we have that $\zb k+\zb h\notin H^{d,0}_R$ for all $\zb k\in H^{d,0}_R$ and $\zb h\in\Lambda(\zb z,M)^\perp\setminus\{\zb 0\}$. 
	Accordingly, we modify the summation index set for $\zb j$ and we estimate the summands
	\begin{align*}
	\theta^2_{\alpha,\beta}(\zb k,\zb z,M)&\lesssim
	\sum_{\zb j\in  \N_0^d\setminus J_R^{d,0}}2^{-2(\alpha\|\zb j\|_1+\beta\|\zb j\|_\infty)}\sum_{\substack{\boldh\in \Lambda(\zb z,M)^{\perp}\\\boldh\neq \boldzero}}\varphi_{\zb j}(\zb k+\zb h).
	\end{align*}
	We apply Lemma \ref{lem:chi_estimate} and get
	\begin{align*}
	\theta^2_{\alpha,\beta}(\zb k,\zb z,M)&\lesssim
	2^{-R}\sum_{\zb j\in  \N_0^d\setminus J_R^{d,0}}2^{-((2\alpha-1)\|\zb j\|_1+\beta\|\zb j\|_\infty)}.
	\end{align*}
	Taking Lemma \ref{lem:summation}  into account, the assertion follows.
\end{proof}

\begin{Lemma}\label{lem:summation}
Let the smoothness parameters $\alpha,\beta\in\R$, $\beta\le0$, $\alpha+\beta>1/2$, and the refinement $R\ge 1$
be given. Then, we estimate
\begin{align*}
 \sum_{\zb j\in\N_0^d\setminus J_R^{d,T}}2^{-((2\alpha-1)\|\zb j\|_1+2\beta\|\zb j\|_\infty)}&\lesssim
\begin{cases}
 2^{-(2\alpha-1+2\beta)R}&:\;T\leq -\frac{\beta}{\alpha} \text{ and }\beta<0,\\
 2^{-(2\alpha-1)R}R^{d-1}&:\;T=\beta=0.
\end{cases}
\end{align*}
\end{Lemma}
\begin{proof}
In the proof of \cite[Theorem 4]{KnDiss} one finds the following estimate
\begin{align*}
 \sum_{\zb j\in\N_0^d\setminus J_R^{d,T}}2^{-t\|\zb j\|_1+s\|\zb j\|_\infty}\lesssim
\begin{cases}
 2^{(s-t)R}&:\;T<\frac{s}{t},\\
 R^{d-1}2^{(s-t+(Tt-s)\frac{d-1}{d-T})R}&:\;T\ge\frac{s}{t}
\end{cases}
\end{align*}
for $s<t$ and $t\ge0$.
Accordingly, we apply this result setting $s:=-2\beta$ and $t:=2\alpha-1$.
We require $\beta\le0$ and obtain the necessity $\alpha+\beta>1/2$ from the conditions $s<t$ and $t\ge0$.
Moreover, we set the parameter $T:=-\beta/\alpha$. This yields
\begin{align*}
 T&=\frac{s}{t+1}
\begin{cases}
 =0&:\;0=s=\beta,\\
 <\frac{s}{t}&:\;0<s=-2\beta.
\end{cases}
\end{align*}
Consequently, we achieve the assertion.
\end{proof}

\begin{Theorem} \label{thm:alias_err:dyadic}
Let the smoothness parameters $\alpha>\frac{1}{2}$, $\beta\leq 0$, $\gamma\ge0$ with $\alpha+\beta>\max\{\gamma,\frac{1}{2}\} $, the dimension $d\in\N$,  $d\ge 2$, and the refinement $R\ge1$, be given.
In addition, we assume that $\Lambda(\zb z,M)$ is a reconstructing rank-1 lattice for $H_R^{d,0}$.
We estimate the error of the sampling operator $\mathrm{Id}-S_{H_R^{d,0}}^{\Lambda(\zb z,M)}$ by
\begin{align*}
M^{-(\alpha+\beta-\gamma)/2}
&\lesssim\|\mathrm{Id}-S_{H_R^{d,0}}^{\Lambda(\zb z,M)}|\mathcal{H}^{\alpha,\beta}(\T^d)\rightarrow \mathcal{H}^{\gamma}_{\mathrm{mix}}(\T^d)\|\lesssim 2^{-(\alpha+\beta-\gamma) R}\begin{cases}
R^{(d-1)/2}&:\;\beta=0,\\
1&:\;\beta<0.
\end{cases}
\intertext{If $\Lambda(\zb z ,M)$ is constructed by the CBC strategy \cite[Tab.\ 3.1]{kaemmererdiss} we continue}
&\lesssim M^{-(\alpha+\beta-\gamma)/2}(\log M)^{\frac{d-2}{2}(\alpha+\beta-\gamma)}\begin{cases}
(\log M)^{(d-1)/2} &:\;\beta =0,\\
1 &:\;\beta<0.
\end{cases}
\end{align*}
\end{Theorem}
\begin{proof}
The lower bound was discussed in Theorem \ref{thm:lower_bound_omega_alpha_beta}.
We apply the triangle inequality and split up the error of the sampling operator 
into the error of the best approximation and the aliasing error.
The error of the projection operator $S_{H_R^{d,0}}$ can be easily estimated using
\begin{align}
\nonumber \|f-S_{H_R^{d,0}}f| \mathcal{H}^\gamma_{\mathrm{mix}}(\T^d)\| = & \Big(\sum_{\boldk \notin H^{d,0}_R}(1+\|\boldk\|_2^2)^{\gamma}|\hat{f}_{\boldk}|^2\Big)^{\frac{1}{2}}\\
\leq & \sup_{k\notin  H^{d,0}_R}\Big(\frac{1}{(1+\|\boldk\|_2^2)^{\beta}\prod_{s=1}^d(1+|k_s|^2)^{\alpha-\gamma}}\Big)^\frac{1}{2} \label{eq:abweighta}\\
&\nonumber \Big(\sum_{\boldk \notin H^{d,0}_R}(1+\|\boldk\|_2^2)^{\beta}\Big[\prod_{s=1}^d(1+|k_s|^2)^{\alpha}\Big]|\hat{f}_{\boldk}|^2\Big)^{\frac{1}{2}}.
\end{align}
It is easy to check that \eqref{eq:abweighta} becomes maximal at the peaks of the hyperbolic cross. Therefore we obtain
\begin{align*}
\|f-S_{H_R^{d,0}}f| \mathcal{H}^\gamma_{\mathrm{mix}}(\T^d)\|\lesssim 2^{-(\alpha+\beta-\gamma)R}\|f|\mathcal{H}^{\alpha,\beta}(\T^d)\|.
\end{align*}

The aliasing error
fulfills
\begin{align}\nonumber
\|S_{H_R^{d,0}}f-S_{H_R^{d,0}}^{\Lambda(\zb z,M)}f|\mathcal{H}^{\gamma}_{\mathrm{mix}}(\T^d)\|^2
=\sum_{\boldk\in H_R^{d,0}}\Big[\prod_{s=1}^d(1+|k_s|^2)^{\gamma}\Big]\Big|\sum_{\substack{\boldh\in \Lambda(\zb z,M)^{\perp}\\\boldh\neq \zb 0}}\hat{f}_{\boldk+\boldh}\Big|^2
\end{align}
Applying H\"older's inequality twice yields
\begin{align}
\|&S_{H_R^{d,0}}f-S_{H_R^{d,0}}^{\Lambda(\zb z,M)}f|\mathcal{H}^{\gamma}_{\mathrm{mix}}(\T^d)\|^2\nonumber\\
&\leq \sum_{\boldk\in H_R^{d,0}}\Big[\prod_{s=1}^d(1+|k_s|^2)^{\gamma}\Big]\Big(\underbrace{\sum_{\substack{\boldh\in \Lambda(\zb z,M)^{\perp}\\\boldh\neq \boldzero}}(1+\|\boldk+\boldh\|_2^2)^{-\beta}\prod_{s=1}^d(1+|k_s+h_s|^2)^{-\alpha}}_{=:\theta^2_{\alpha,\beta}(\zb k,\zb z,M),\text{ cf. \eqref{def:theta}}}\Big).\nonumber\\
&\Big(\sum_{\substack{\boldh\in \Lambda(\zb z,M)^{\perp}\\\boldh\neq \boldzero}}(1+\|\boldk+\boldh\|_2^2)^{\beta}\Big[\prod_{s=1}^d(1+|k_s+h_s|^2)^{\alpha}\Big]|\hat{f}_{\boldk+\boldh}|^2\Big)\nonumber\\
&\leq \sup_{\boldk\in H_R^{d,0}}\Big[\prod_{s=1}^d(1+|k_s|^2)^{\gamma}\Big]\theta^2_{\alpha,\beta}(\zb k,\zb z,M)\nonumber\\
&\Big(\sum_{\boldk\in H_R^{d,0}}\sum_{\substack{\boldh\in \Lambda(\zb z,M)^{\perp}\\\boldh\neq \boldzero}}(1+\|\boldk+\boldh\|_2^2)^{\beta}\Big[\prod_{s=1}^d(1+|k_s+h_s|^2)^{\alpha}\Big]|\hat{f}_{\boldk+\boldh}|^2\Big)\nonumber\\
&\leq \sup_{\boldh\in H_R^{d,0}}\Big[\prod_{s=1}^d(1+|h_s|^2)^{\gamma}\Big]\sup_{\boldk\in H_R^{d,0}}\theta^2_{\alpha,\beta}(\zb k,\zb z,M)\|f|\mathcal{H}^{\alpha,\beta}(\T^d)\|^2\label{eq:47}
\end{align}
since $\Lambda(\zb z,M)$ is a reconstructing rank-1 lattice for $H_R^{d,0}$ and, consequently,
the sets
$\{\zb k+\zb h\in\Z^d\colon \zb h\in\Lambda(\zb z,M)^\perp\}$, $\zb k\in H_R^{d,0}$,
do not intersect.
We apply Lemma \ref{lem:theta} and take the upper bound
$$
\sup_{\zb k\in H_R^{d,0}}\prod_{s=1}^d(1+|k_s|^2)^\gamma\lesssim\sup_{\zb j\in J_R^{d,0}}2^{2\gamma\|\zb j\|_1}\lesssim 2^{2\gamma R}
$$
into account. We achieve
\begin{align*}
\|S_{H_R^{d,0}}f-S_{H_R^{d,0}}^{\Lambda(\zb z,M)}f|\mathcal{H}^{\gamma}_{\mathrm{mix}}(\T^d)\|
&\lesssim \|f|\mathcal{H}^{\alpha,\beta}(\T^d)\|\,2^{-(\alpha+\beta-\gamma) R} \begin{cases}
R^{\frac{d-1}{2}}&:\;\beta=0,\\
1&:\;\beta<0
\end{cases}
\end{align*}
and, in conjunction with Lemma \ref{lem:universal_bounds_M}, the second assertion of the theorem.
\end{proof}
\begin{Remark}
The basic improvement in the error analysis compared to \cite{KaPoVo14} is provided by applying Lemma \ref{lem:chi_estimate} in \eqref{eq:47}. Here, the information about the cardinality of the dual lattice intersected with rectangular boxes yields sharp main rates coinciding with the lower bounds given in Theorem \ref{thm:lower_bound_omega_alpha_beta}. From that viewpoint this technique improves also the asymptotical main rates obtained in \cite{KuSlWo06} for the $L_2(\T^d)$ approximation error. In case $\beta<0$ and $\gamma=0$ the result above behaves not optimal compared to the result obtained in \cite{KaPoVo14} where a Korobov type lattice is used. The authors there obtain no logarithmic dependence in $M$. The main reason for that issue is the probably technical limitation in Lemma \ref{lem:chi_estimate} discussed in Remark \ref{rem:countingissue} that does not allow us to use energy-type hyperbolic crosses as index sets, here.
\end{Remark}

\noindent Due to the embedding $\mathcal{H}^{\gamma}_{\mathrm{mix}}(\T^d)\hookrightarrow \mathcal{H}^{\gamma}(\T^d)$ we obtain the following proposition.

\begin{Proposition} \label{prop:err:dyadic}
Let the smoothness parameters $\alpha>\frac{1}{2}$, $\beta\leq 0$, $\gamma\geq0$ with $\alpha+\beta>\max\{\gamma,\frac{1}{2}\} $, the dimension $d\in\N$,  $d\ge 2$, and the refinement $R\ge 1$, be given.
In addition, we assume that $\Lambda(\zb z,M)$ is a reconstructing rank-1 lattice for $H_R^{d,0}$
constructed by the CBC strategy \cite[Tab.\ 3.1]{kaemmererdiss}.
We estimate the error of the sampling operator $\mathrm{Id}-S_{H_R^{d,0}}^{\Lambda(\zb z,M)}$ by
\begin{align*}
M^{-(\alpha+\beta-\gamma)/2}
&\lesssim\|\mathrm{Id}-S_{H_R^{d,0}}^{\Lambda(\zb z,M)}|\mathcal{H}^{\alpha,\beta}(\T^d)\rightarrow \mathcal{H}^{\gamma}(\T^d)\|\lesssim 2^{-(\alpha+\beta-\gamma) R}\begin{cases}
R^{(d-1)/2}&:\;\beta=0,\\
1&:\;\beta<0
\end{cases}
\\ &\lesssim M^{-(\alpha+\beta-\gamma)/2}(\log M)^{(d-2)(\alpha+\beta-\gamma)/2}\begin{cases}
(\log M)^{(d-1)/2}&:\;\beta =0,\\
1&:\;\beta<0.
\end{cases}
\end{align*}
\par \vspace{-1.7\baselineskip}
\qed
\end{Proposition}

For $2<q<\infty$ the embedding $$\mathcal{H}^{\frac{1}{2}-\frac{1}{q}}(\T^d)\hookrightarrow L_q(\T^d)$$ (see \cite{MR891189}, 2.4.1) extends the last theorem to target spaces $L_q(\T^d)$.

\begin{Proposition} \label{prop:err:dyadic:Lq}
Let the smoothness parameters $\alpha>\frac{1}{2}$ and $\beta\leq 0$ with $\alpha+\beta>\frac{1}{2}$, $2<q<\infty$. Let the dimension $d\in\N$,  $d\ge 2$, and the refinement $R\ge1$, be given.
In addition, we assume that $\Lambda(\zb z,M)$ is a reconstructing rank-1 lattice for $H_R^{d,0}$
constructed by the CBC strategy \cite[Tab.\ 3.1]{kaemmererdiss}.
We estimate the error of the sampling operator $\mathrm{Id}-S_{H_R^{d,0}}^{\Lambda(\zb z,M)}$ by
\begin{align*}
\|\mathrm{Id}&-S_{H_R^{d,0}}^{\Lambda(\zb z,M)}|\mathcal{H}^{\alpha,\beta}(\T^d)\rightarrow L_q(\T^d)\|\lesssim 2^{-(\alpha+\beta-(\frac{1}{2}-\frac{1}{q})) R}\begin{cases}
R^{(d-1)/2}&:\;\beta=0,\\
1&:\;\beta<0
\end{cases}\\
&\lesssim M^{-(\alpha+\beta-(\frac{1}{2}-\frac{1}{q}))/2}(\log M)^{\frac{d-2}{2}(\alpha+\beta-(\frac{1}{2}-\frac{1}{q}))} \begin{cases}
(\log M)^{(d-1)/2}&:\;\beta=0,\\
1&:\;\beta<0.
\end{cases}
\end{align*}
\par \vspace{-1.7\baselineskip}
\qed
\end{Proposition}

In addition to $L_q(\T^d)$, $2<q<\infty$, we study the case $q=\infty$.
For technical reasons we estimate the sampling error with respect to the $d$-dimensional Wiener algebra
$$\mathcal{A}(\T^d):=\{f\in L_1(\T^d)\colon\sum_{\boldk\in\Z^d}|\hat{f}_{\boldk}|<\infty\}$$
and subsequently we use the embedding  $\mathcal{A}(\T^d)\hookrightarrow  C(\T^d)\hookrightarrow L_{\infty}(\T^d)$.
\begin{Theorem} \label{thm:alias_infty_err:dyadic}
Let the smoothness parameters $\alpha>\frac{1}{2}$ and $\beta \leq 0$ with $\alpha+\beta>\frac{1}{2}$, the dimension $d\in\N$,  $d\ge 2$, and the refinement $R\in\R$, $R\ge1$, be given.
In addition, we assume that $\Lambda(\zb z,M)$ is a reconstructing rank-1 lattice for $H_R^{d,T}$ with $T:=-\frac{\beta}{\alpha}$
constructed by the CBC strategy \cite[Tab.\ 3.1]{kaemmererdiss}.
We estimate the error of the sampling operator $\mathrm{Id}-S_{H_R^{d,T}}^{\Lambda(\zb z,M)}$ by
\begin{align*}
\|\mathrm{Id}-S_{H_R^{d,T}}^{\Lambda(\zb z,M)}|\mathcal{H}^{\alpha,\beta}(\T^d)\rightarrow\mathcal{A}(\T^d)\|
&\lesssim 2^{-(\alpha+\beta-\frac{1}{2})R}\begin{cases}
R^{\frac{d-1}{2}} &:\;\beta =0,\\
1 &:\;\beta<0
\end{cases}
\\ &\lesssim M^{-(\alpha+\beta-\frac{1}{2})/2}\begin{cases}
(\log M)^{\frac{d-2}{2}(\alpha-\frac{1}{2})+\frac{d-1}{2}} &:\;\beta =0,\\
1 &:\;\beta<0.
\end{cases}
\end{align*}
\end{Theorem}
\begin{proof}
Again we use the triangle inequality and split up the error of the sampling operator 
into the error of the truncation error and the aliasing error. The truncation error fulfills
\begin{align}
\|f-S_{H_R^{d,T}}f|\mathcal{A}(\T^d)\| \lesssim \|f|\mathcal{H}^{\alpha,\beta}(\T^d)\| 2^{-(\alpha+\beta-\frac{1}{2})R}\begin{cases}
R^{\frac{d-1}{2}} &:\;\beta=0,\\
1&:\;\beta<0.
\end{cases}\label{eq:trunc_err_wieneralg}
\end{align}
 For completeness we give a short proof. Applying the orthogonal projection property of $S_{H_R^{d,T}}f$ we obtain 
\begin{align*}
\|&f-S_{H_R^{d,T}}f|\mathcal{A}(\T^d)\|=\sum_{\boldk\notin H_R^{d,T}}|\hat{f}_\boldk|\\
&\leq\Big(\sum_{\boldk\notin H_R^{d,T}}(1+\|\boldk\|_2^2)^{-\beta}\prod_{s=1}^d(1+|k_s|^2)^{-\alpha}\Big)^{\frac{1}{2}} \Big(\sum_{\boldk\notin H_R^{d,T}}(1+\|\boldk\|_2^2)^{\beta}\Big[\prod_{s=1}^d(1+|k_s|^2)^{\alpha}\Big]|\hat{f}_\boldk|^2\Big)^{\frac{1}{2}}.
\end{align*}
Decomposing the first sum into dyadic blocks yields
\begin{align}
\|f-S_{H_R^{d,T}}f|\mathcal{A}(\T^d)\|&\leq \Big(\sum_{\boldj\notin J_R^{d,T}}\sum_{\boldk\in Q_\boldj}(1+\|\boldk\|_2^2)^{-\beta}\prod_{s=1}^d (1+|k_s|^2)^{-\alpha}\Big)^{\frac{1}{2}}\|f|\mathcal{H}^{\alpha,\beta}(\T^d)\|\label{eq:a_error_eq1}\\
&\lesssim \Big(\sum_{\boldj\notin J_R^{d,T}}2^{-2\alpha\|\boldj\|_1-2\beta\|\boldj\|_{\infty}}\sum_{\boldk\in Q_\boldj}1\Big)^{\frac{1}{2}}\|f|\mathcal{H}^{\alpha,\beta}(\T^d)\|\nonumber\\
&\lesssim \Big(\sum_{\boldj\notin J_R^{d,T}}2^{-(2\alpha-1)\|\boldj\|_1-2\beta\|\boldj\|_{\infty}}\Big)^{\frac{1}{2}}\|f|\mathcal{H}^{\alpha,\beta}(\T^d)\|.\nonumber
\end{align}
Applying Lemma \ref{lem:summation} we obtain \eqref{eq:trunc_err_wieneralg}.
The aliasing error behaves as follows
\begin{align*}
\|S_{H_R^{d,T}}f-S_{H_R^{d,T}}^{\Lambda(\zb z,M)}f|\mathcal{A}(\T^d)\|
=\sum_{\boldk\in H_R^{d,T}}\Big|\sum_{\substack{\boldh\in \Lambda(\zb z,M)^{\perp}\\\boldh\neq \boldzero}}\hat{f}_{\boldk+\boldh}\Big|.
\end{align*}
Applying H\"older's inequality twice yields
\begin{align*}
\|&S_{H_R^{d,T}}f-S_{H_R^{d,T}}^{\Lambda(\zb z,M)}f|\mathcal{A}(\T^d)\|\\
&\leq \Big(\sum_{\boldk\in H_R^{d,T}}\sum_{\substack{\boldh\in \Lambda(\zb z,M)^{\perp}\\\boldh\neq \boldzero}}(1+\|\boldk+\boldh\|_2^2)^{-\beta}\prod_{s=1}^d (1+|k_s+h_s|^2)^{-\alpha}\Big)^{\frac{1}{2}}\\
&\Big(\sum_{\boldk\in H_R^{d,T}}\sum_{\substack{\boldh\in \Lambda(\zb z,M)^{\perp}\\\boldh\neq \boldzero}}(1+\|\boldk+\boldh\|_2^2)^{\beta}\prod_{s=1}^d (1+|k_s+h_s|^2)^{\alpha}|\hat{f}_{\boldk+\boldh}|^2\Big)^{\frac{1}{2}}.
\end{align*}
Since $\Lambda(\zb z,M)$ is a reconstructing rank-1 lattice for $H_R^{d,T}$ and, consequently,
the sets\linebreak
$\{\zb k+\zb h\in\Z^d\colon \zb h\in\Lambda(\zb z,M)^\perp\}$, $\zb k\in H_R^{d,T}$,
do not intersect, we obtain
\begin{align*}
\|&S_{H_R^{d,T}}f-S_{H_R^{d,T}}^{\Lambda(\zb z,M)}f|\mathcal{A}(\T^d)\|\\
&\leq \Big(\sum_{\boldk\notin H_R^{d,T}}(1+\|\boldk\|_2^2)^{-\beta}\prod_{s=1}^d (1+|k_s|^2)^{-\alpha}\Big)^{\frac{1}{2}}\Big(\sum_{\boldk\notin H_R^{d,T}}(1+\|\boldk\|_2^2)^{\beta}\prod_{s=1}^d (1+|k_s|^2)^{\alpha}|\hat{f}_{\boldk}|^2\Big)^{\frac{1}{2}}\\
&\leq \Big(\sum_{\boldk\notin H_R^{d,T}}(1+\|\boldk\|_2^2)^{-\beta}\prod_{s=1}^d (1+|k_s|^2)^{-\alpha}\Big)^{\frac{1}{2}}\|f|\mathcal{H}^{\alpha,\beta}(\T^d)\|.
\end{align*}
Now we are in the same situation as in \eqref{eq:a_error_eq1}. Therefore we achieve
\begin{align*}
\|S_{H_R^{d,T}}f-S_{H_R^{d,T}}^{\Lambda(\zb z,M)}f|\mathcal{A}(\T^d)\|&\lesssim \|f|\mathcal{H}^{\alpha,\beta}(\T^d)\| 2^{-(\alpha+\beta-\frac{1}{2})R}\begin{cases}
R^{\frac{d-1}{2}}&:\;\beta=0,\\
1&:\;\beta<0.
\end{cases}
\end{align*}
Here, we would like to particularly mention that the aliasing error has the same order as the truncation error.
\end{proof}
\begin{Proposition} \label{thm:alias_err_infty:dyadic}
Let the smoothness parameter $\alpha>\frac{1}{2}$ and $\beta \leq 0$ with $\alpha+\beta>\frac{1}{2}$, the dimension $d\in\N$,  $d\ge 2$, and the refinement $R\ge1$, be given.
In addition, we assume that $\Lambda(\zb z,M)$ is a reconstructing rank-1 lattice for $H_R^{d,T}$ with $T:=-\frac{\beta}{\alpha}$
constructed by the CBC strategy \cite[Tab.\ 3.1]{kaemmererdiss}.
We estimate the error of the sampling operator $\mathrm{Id}-S_{H_R^{d,T}}^{\Lambda(\zb z,M)}$ by
\begin{align*}
\|\mathrm{Id}-S_{H_R^{d,T}}^{\Lambda(\zb z,M)}|\mathcal{H}^{\alpha,\beta}(\T^d)\rightarrow L_{\infty}(\T^d)\|
&\lesssim 2^{-(\alpha+\beta-\frac{1}{2})R}\begin{cases}
R^{\frac{d-1}{2}} &:\;\beta =0,\\
1 &:\;\beta<0
\end{cases}
\\ &\lesssim M^{-(\alpha+\beta-\frac{1}{2})/2}\begin{cases}
(\log M)^{\frac{d-2}{2}(\alpha-\frac{1}{2})+\frac{d-1}{2}} &:\;\beta =0,\\
1 &:\;\beta<0.
\end{cases}
\end{align*}
\end{Proposition}
\begin{Remark}
In case $\beta<0$ the technique used in the proof of Theorem \ref{thm:alias_infty_err:dyadic} and Proposition \ref{thm:alias_err_infty:dyadic} allows it to benefit from smaller index sets $H^{d,T}_R$ with $T>0$, so called energy-type hyperbolic crosses. Therefore, we obtain no logarithmic dependencies in the error rate.
\end{Remark}
\section{The two-dimensional case}\label{sec:2d}
In this chapter we restrict our considerations to two-dimensional approximation problems, i.e., the dimension $d=2$ is fixed.
We collect some basic facts from above on this special case.
\begin{Lemma}\label{lem:2dimlattice}
Let $R\ge 0$, and $T\in[0,1)$ be given.
Each reconstructing rank-1 lattice $\Lambda(\zb z,M)$ for the frequency index set $H_{R}^{2,T}\subset\Z^2$ fulfills
\begin{itemize}
\item $M\ge 2^{2\floor{R}}$,
\item $\Lambda(\zb z,M)$ is a reconstructing rank-1 lattice for the tensor product grid \newline
$G^2_R:=(-2^{\floor{R}-1},2^{\floor{R}-1}]^2\cap\Z^2$.
\end{itemize}
Moreover, there exist reconstructing rank-1 lattices $\Lambda(\zb z,M)$ for the frequency index sets $H_{R}^{2,T}$ that
fulfill $M=(1+3\cdot2^{\ceil{R}-1})2^{\ceil{R}}\le 2^{2R+3}$.
\end{Lemma}
\begin{proof}
The proof follows from \cite[Theorem 3.5 and Lemma 3.7]{KaKuPo10} and the embeddings $H_R^{2,T}\subset H_R^{2,0}$ for $T\ge0$, which is direct consequence of the definition.
\end{proof}

We interpret the last lemma. The reconstruction property of reconstructing rank-1 lattices $\Lambda(\zb z,M)$ for two-dimensional hyperbolic crosses $H_R^{2,T}\subset (-2^{R},2^{R}]^2\cap\Z^2$
implies automatically that the rank-1 lattices $\Lambda(\zb z,M)$ are reconstructing rank-1 lattices for only mildly lower expanded full grids $(-2^{\floor{R}-1},2^{\floor{R}-1}]^2\cap\Z^2$.
Accordingly, in the sense of sampling numbers it seems appropriate to use a rank-1 lattice sampling in combination with tensor product grids as frequency index sets in order
to even approximate functions of dominating mixed smoothness in dimensions $d=2$.
Thus, we consider the sampling operator $S^{\Lambda(\zb z,M)}_{G_R^2}$, cf. \eqref{eqn:sampling_operator}.

\begin{Lemma}\label{lem:summation_linfty_ball}
Let $a\in \R$, $0<a<1$ and $L\in\N$ be given. Then we estimate
$$\sum_{\substack{\zb j\in\N_0^2\\\|\zb j\|_\infty\ge L}}a^{\|\zb j\|_1}\le \frac{2-a^L}{(1-a)^2}a^L\le C_a\cdot a^L.$$
\end{Lemma}
\begin{proof}
We evaluate the geometric series and get
\begin{align*}
\sum_{\substack{\zb j\in\N_0^2\\\|\zb j\|_\infty\ge L}}a^{\|\zb j\|_1}
&=\sum_{j_1=0}^{L-1}a^{j_1}\sum_{j_2=L}^{\infty}a^{j_2}+\sum_{j_2=0}^{L-1}a^{j_2}\sum_{j_1=L}^{\infty}a^{j_1}+\sum_{j_1=L}^{\infty}a^{j_1}\sum_{j_2=L}^{\infty}a^{j_2}
\\
&=\left(\frac{1-a^L}{1-a}+\frac{1-a^L}{1-a}+\frac{a^L}{1-a}\right)\frac{a^L}{1-a}.
\end{align*}
\end{proof}
\begin{Theorem} \label{thm:alias_2d_err:dyadic}
Let the smoothness parameter $\alpha>\frac{1}{2},\gamma \geq  0$ with $\alpha>\gamma$ and the refinement $R\ge0$, be given.
In addition, we assume that $\Lambda(\zb z,M)$ is a reconstructing rank-1 lattice for $G_R^2$ with $M\asymp2^{2R}$.
We estimate the error of the sampling operator $\mathrm{Id}-S^{\Lambda(\zb z,M)}_{G_R^2}$ by
\begin{align*}
\|\mathrm{Id}-S^{\Lambda(\zb z,M)}_{G_R^2}|\mathcal{H}^{\alpha}_{\mathrm{mix}}(\T^2)\rightarrow \mathcal{H}^{\gamma}(\T^2)\|
\asymp M^{-(\alpha-\gamma)/2}.
\end{align*}

\end{Theorem}
\begin{proof}
The lower bound goes back to Theorem \ref{thm:lower_bound_omega_alpha_beta}.
The proof of the upper bound is similar to the proof of Theorem \ref{thm:alias_err:dyadic}. The main difference is that we use the full grid $G_R^2$ instead of $H^{2,0}_R$ here. This yields for the projection
$$\|\mathrm{Id}-S_{G_R^2}|\mathcal{H}^{\alpha}_{\mathrm{mix}}(\T^2)\rightarrow \mathcal{H}^{\gamma}(\T^2)\| \lesssim  M^{-(\alpha-\gamma)/2}.$$
The estimation for the aliasing error $\|S_{G_R^2}f-S^{\Lambda(\zb z,M)}_{G_R^2}f|\mathcal{H}^{\gamma}(\T^2)\|$ is also very similar to \eqref{thm:alias_err:dyadic}. We follow the proof line by line with the mentioned modification and come to the estimation
\begin{align*}
\|S_{G_R^2}f&-S^{\Lambda(\zb z,M)}_{G_R^2}f|\mathcal{H}^{\gamma}(\T^2)\|\\
&\leq 
\sup_{\zb k\in G_R^2}\Big((1+\|\zb k\|_2^2)^{\gamma}\sum_{\zb j\in\N_0^d}\sum_{\substack{\zb h\in \Lambda(\zb z,M)^{\perp}\\\zb h\neq \zb 0}}\varphi_{\zb j}(\zb k+\zb h)\prod_{i=1}^d(1+|k_i+ h_i|^2)^{-\alpha}\Big)^{\frac{1}{2}}\|f|\mathcal{H}^{\alpha}_{\mathrm{mix}}(\T^2)\|.
\end{align*}
Due to the reproduction property for $G_R^2$ the sum over $\zb j$ breaks down to
\begin{align*}
\|S_{G_R^2}f&-S^{\Lambda(\zb z,M)}_{G_R^2}f|\mathcal{H}^{\gamma}(\T^2)\|\\
&\lesssim \sup_{\zb k\in G_R^2}\Big((1+\|\zb k\|_2^2)^{\gamma}\sum_{\|\zb j\|_{\infty}>\floor{R}}2^{-2\alpha\|\zb j\|_1}\sum_{\substack{\zb h\in \Lambda(\zb z,M)^{\perp}\\\zb h\neq \zb 0}}\varphi_{\zb j}(\zb k+\zb h)\Big)^{\frac{1}{2}}\|f|\mathcal{H}_{\mathrm{mix}}^{\alpha}(\T^2)\|.
\end{align*}
Next, we recognize
\begin{equation}
\sup_{\zb k\in G_R^2}(1+\|\zb k\|^2_2)^{\frac{\gamma}{2}}\lesssim 2^{\gamma R}.\label{eq:isoweightestimation} 
\end{equation}

Using $H_{R-2}^{d,0}\subset G_R^2$, we obtain $\Lambda(\zb z,M)^\perp\cap H_{R-2}^{d,0}=\{\zb 0\}$. We apply Lemma \ref{lem:chi_estimate} and employ $R-1\le\floor{R}\le\|\zb j\|_\infty\le\|\zb j\|_1$ to see
\begin{align*}
\|S_{G_R^2}f&-S^{\Lambda(\zb z,M)}_{G_R^2}f|\mathcal{H}^{\gamma}(\T^2)\|\\
&\lesssim 2^{\gamma R}\Big(2^{-R} \sum_{\|\zb j\|_{\infty}>\floor{R}}2^{-(2\alpha-1) \|\zb j\|_1}\Big)^{\frac{1}{2}}\|f|\mathcal{H}^{\alpha}_{\mathrm{mix}}(\T^2)\|.\\
\intertext{Applying Lemma
 \ref{lem:summation_linfty_ball} yields}
\|S_{G_R^2}f-S^{\Lambda(\zb z,M)}_{G_R^2}f|\mathcal{H}^{\gamma}(\T^2)\|&\lesssim 2^{-(\alpha-\gamma)R}\|f|\mathcal{H}^{\alpha}_{\mathrm{mix}}(\T^2)\|\\
&\lesssim M^{-(\alpha-\gamma)/2}\|f|\mathcal{H}^{\alpha}_{\mathrm{mix}}(\T^2)\|.
\end{align*}
\end{proof}
\begin{Remark}\label{rem:err_2d_Hmix}
This method does not work for $\mathcal{H}^{\gamma}_{\mathrm{mix}}(\T^2)$ as target space. Here the estimation of the mixed weight, similar to \eqref{eq:isoweightestimation} implies a worse main rate for the asymptotic behavior of $\|S_{G_R^2}f-S_{G_R^2}^{\Lambda(\zb z,M)}f|\mathcal{H}^{\gamma}_{\mathrm{mix}}(\T^2)\|$. Here we have to use $H^{2,0}_R$ as index set for our trigonometric polynomials and therefore Theorem \ref{thm:alias_err:dyadic} is the best we have in this situation.
\end{Remark}
\begin{Theorem}
Let the smoothness parameter $\alpha>\frac{1}{2}$ and the refinement $R\ge0$ be given.
In addition, we assume that $\Lambda(\zb z,M)$ is a reconstructing rank-1 lattice for $G_R^2$ with $M\asymp2^{2R}$.
We estimate the error of the sampling operator $\mathrm{Id}-S_{G_R^2}^{\Lambda(\zb z,M)}$ by
\begin{align*}
\|\mathrm{Id}-S_{G_R^2}^{\Lambda(\zb z,M)}|\mathcal{H}^{\alpha}_{\mathrm{mix}}(\T^2)\rightarrow\mathcal{A}(\T^2)\|
&\lesssim M^{-(\alpha-\frac{1}{2})/2}.
\end{align*}
\end{Theorem}
\begin{proof}
The result is a consequence of replacing $H^{2,0}_R$ by $G_R^2$ in the proof of Theorem \ref{thm:alias_infty_err:dyadic}.
\end{proof}
\begin{Proposition}\label{prop:err_2d_Linf}
Let the smoothness parameter $\alpha>\frac{1}{2}$ and the refinement $R\ge0$, be given.
In addition, we assume that $\Lambda(\zb z,M)$ is a reconstructing rank-1 lattice for $G_R^2$ with $M\asymp2^{2R}$.
We estimate the error of the sampling operator $\mathrm{Id}-S_{G_R^2}^{\Lambda(\zb z,M)}$ by
\begin{align*}
\|\mathrm{Id}-S_{G_R^2}^{\Lambda(\zb z,M)}|\mathcal{H}^{\alpha}_{\mathrm{mix}}(\T^2)\rightarrow L_{\infty}(\T^2)\|
&\lesssim M^{-(\alpha-\frac{1}{2})/2}.
\end{align*}
\par \vspace{-1.7\baselineskip}
\qed
\end{Proposition}
\begin{figure}[tb]
\begin{center}
\begin{tikzpicture}[scale=0.7]
\begin{axis}[
axis x line=none,
axis y line=none,
ticks=none,
xmin=-5,
xmax=5,
ymin=-5,
ymax=5,
height=10cm,
width=10cm,
legend pos=south east,
xticklabel style={/pgf/number format/.cd,fixed}
]

\addplot[mark=none, color=black, domain=0.1:5,thick]
plot[samples=200,smooth] {4/x};
\addplot[mark=none, color=black, domain=0.1:5,thick]
plot[samples=200,smooth] {-4/x};
\addplot[mark=none, color=black, domain=-5:-0.1,thick]
plot[samples=200,smooth] {-4/x};
\addplot[mark=none, color=black, domain=-5:-0.1,thick]
plot[samples=200,smooth] {4/x};
\fill[color=gray] (axis cs:-2,-2) -- (axis cs:2,-2) -- (axis cs:2,2) -- (axis cs:-2,2) -- (axis cs:-2,-2);
\fill[color=lightgray] (axis cs:-1,-1) -- (axis cs:1,-1) -- (axis cs:1,1) -- (axis cs:-1,1) -- (axis cs:-1,-1);

\node at (axis cs:0,0) {$ B_n$};
\node at (axis cs:0,1.4) {$\mathcal{D}(B_n)$};
\end{axis}

\end{tikzpicture} 
\end{center}
\caption{Relations between $B_n,\mathcal{D}(B_n)$ and a hyperbolic cross of size $\delta b_n$.}\label{fig:fibonacci_full_in_hc}
\end{figure}
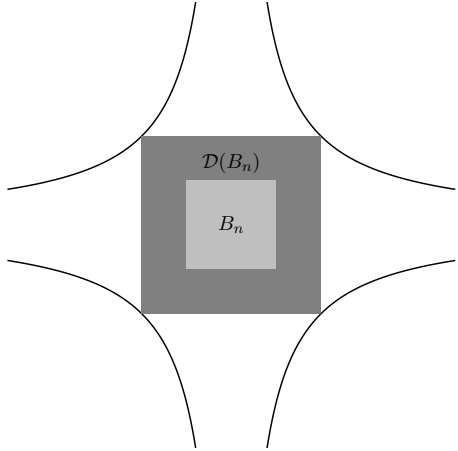
Now we come to the second very special property of the 2-dimensional situation. Here we know closed formulas for lattices that are reconstructing for $H^{2,0}_R$ (and $G_R^2$). The well studied Fibonacci lattice $F_n=\Lambda(\zb z,b_n)$,
where $\boldz=(1,b_{n-1})$ and $M=b_n$ gives a universal reconstructing rank-1 lattice for index sets considered in this chapter.
The Fibonacci numbers $b_n$ are defined iteratively by
$$b_0=b_1=1,\quad b_n=b_{n-1}+b_{n-2},\;n\geq 2.$$ 
Since the size of the Fibonacci lattice depends on $M=b_n$, we go the other way around. For a fixed refinement $n\in
\N$ we choose a suitable rectangle $B_n$ for which the reproduction property \eqref{eqn:dual_lattice_and_difference_set}
is fulfilled. Let us start with the box
$$B_n:=\left[-\left\lfloor C \sqrt{b_n}\right\rfloor,\left\lfloor C\sqrt{b_n}\right\rfloor\right]^2\cap\Z^2,$$
where $C>0$ is a suitable constant.
Obviously, the difference set of such a box fulfills
$$\mathcal{D}(B_n)= \left[-2\left\lfloor C \sqrt{b_n}\right\rfloor,2\left\lfloor C\sqrt{b_n}\right\rfloor\right]^2\cap\Z^2.$$
It is known (see Lemma IV.2.1 in \cite{Tem93}), that there is a $\delta>0$ such that for all frequencies of the dual lattice $F_n^{\perp}$ of $F_n$ 
$$\prod_{s=1}^2 \max\{1,|h_s|\}\geq \delta b_n$$
holds.
For that reason we find a $C>0$ (depending only on $\delta$) such that the property
$$\mathcal{D}(B_{n})\cap F_n^{\perp} =\{0\}$$
is fulfilled for all $n\in \N$ (see Figure \ref{fig:fibonacci_full_in_hc}), which guarantees the reproduction property for the index set $B_n$.
Additionally we have
$|B_n| \asymp b_n.$
Therefore, the Fibonacci lattice fulfills the properties mentioned in Lemma \ref{lem:2dimlattice}.

\section{Further comments}
\subsection{Minkowski's theorem in Section \ref{sec:lower_bounds}}
\begin{Remark}\label{rem:minkowski}
In order to show the lower bounds in Theorem~\ref{thm:lower_bound_omega_alpha_beta},
one may alternatively use Minkowski's theorem
instead of the construction in Lemma \ref{lem:aliasing_axis_cross}.
Then the main rate in $M$ is identical but one obtains an additional
factor that decreases exponentially in the dimension $d$ in the lower bound.
\end{Remark}

\subsection{Hyperbolic cross property in Section \ref{sec:upper_bounds}}
The following remark is hypothetical since it is an open question whether a lattice with the so-called ``hyperbolic cross property'' exists in $d>2$, cf. Lemma~\ref{lem:universal_bounds_M}.

\begin{Remark}\label{rem:hyp_cross_prop}
Let $\Lambda(\zb z,M)$ be a lattice such that $\Lambda(\zb z,M)^{\perp}\cap H^{d,0}_{2R}=\{0\} $ with $M\asymp 2^{2R}$ holds. We call this property ``hyperbolic cross property". Then
\begin{align*}
\|f-S^{\Lambda(\zb z,M)}_{H^{d,0}_{R}}f|\mathcal{H}^{\gamma}(\T^d)\|&\lesssim
2^{-(\alpha-\gamma)R}R^{\frac{d-1}{2}}\|f|\mathcal{H}^{\alpha}_{\mathrm{mix}}(\T^d)\|\\
&\asymp  M^{-\frac{\alpha-\gamma}{2}}(\log M)^{\frac{d-1}{2}}.
\end{align*}
\end{Remark}
\begin{proof}
Computing the truncation error is straight-forward. For the aliasing error we get
\begin{align*}
\|&S_{H^{d,0}_{R}}f-S^{\Lambda(\zb z,M)}_{H^{d,0}_{R}}f|\mathcal{H}^{\gamma}(\T^d)\|\\
\leq& \sup_{\boldk\in H_{R}^{d,0}}\Big((1+\|\boldk\|^2_2)^{\gamma}\sum_{\boldj\notin J^{d,0}_{R}}\sum_{\substack{\zb h\in \Lambda(\zb z,M)^{\perp}\\\zb h\neq \zb 0}}
\varphi_\boldj(\boldk+\boldh)\prod_{s=1}^d(1+|k_s+h_s|^2)^{-\alpha}\Big)^{\frac{1}{2}}\|f|\mathcal{H}^{\alpha}_{
\mathrm{mix}}(\T^d)\|.
\end{align*}
Now we use the fact that the difference set  $\mathcal D(H^{d,0}_{R})$ is contained in $H^{d,0}_{c+2R}$ and therefore,
$\Lambda(\zb z,M)$ is reproducing for $H^{d,0}_{R}$ (the dual lattice is located outside of the difference set). With
the usual calculation we get then
\begin{align*}
\|&S_{H^{d,0}_{R}}f-S^{\Lambda(\zb z,M)}_{H^{d,0}_{R}}f|\mathcal{H}^{\gamma}(\T^d)\|\\
\lesssim& \sup_{\boldk\in
H_{R}^{d,0}}(1+\|\boldk\|^2_2)^{\frac{\gamma}{2}}\Big(\sum_{R<\|\boldj\|_1<2R}2^{-2\alpha\|\bold j\|_1}+\sum_{
\|\boldj\|_1>2R}2^{-2\alpha\|\boldj\|_1}\frac{2^{\|\boldj\|_1}}{2^{2R}}\Big)^{ \frac{1}{2}}\\
\lesssim& 2^{-(\alpha-\gamma)R}R^{\frac{d-1}{2}}\|f|\mathcal{H}^{\alpha}_{\mathrm{mix}}(\T^d)\|\\
\asymp& M^{-\frac{\alpha-\gamma}{2}}(\log M)^{\frac{d-1}{2}}\|f|\mathcal{H}^{\alpha}_{\mathrm{mix}}(\T^d)\|.
\end{align*}
\end{proof}
Unfortunately, if $d>2$ such a lattice is not known. We see that even in this ``ideal'' case we do not get
rid of the $(\log M)^{\frac{d-1}{2}}$. If $d=2$ we get rid of both $\log$s, see
Section \ref{sec:2d}. One reason is that e.g. the Fibonacci lattice has a ``hyperbolic cross property" (cf. Remark
\ref{rem:hyp_cross_prop}). The other reason is that due to the ``half rate'' we can truncate
from a larger set than the hyperbolic cross. In that sense $d=2$ is a very specific case.

\subsection{Energy-norm setting in Section \ref{sec:upper_bounds}}
\begin{Remark}\label{rem:countingissue}
Additionally to the considerations in
Proposition \ref{prop:err:dyadic} it seems natural to treat the cases $\gamma>\beta>0$. One would expect from the theory of sparse
grids that a modification of the hyperbolic cross index sets $H^{d,0}_R$ to energy-norm based hyperbolic crosses
$H^{d,T}_R$ with $T = \frac{\gamma-\beta}{\alpha}$ or a little perturbation of it would
help to reduce logarithmic dependence on $M$. Unfortunately, we are currently not able to improve or even get
equivalent results for that. One reason is that we have no improved results fitting
$H^{d,T}_R$ in Lemma \ref{lem:chi_estimate}. The other reason is that in case $\gamma>0 $ we have not yet found a way
to exploit smoothness that come from the target space such that one can use smaller
index sets than $H^{d,0}_R$ in the error sum. Our standard estimation yields a worse main rate for that.
\end{Remark}

\subsection{Sampling along multiple rank-1 lattices}
Similar to sampling along sparse grids, which are unions of anisotropic full grids,
one may use the union of several rank-1 lattices as sampling set, cf. \cite{Ka15}.
In contrast to the CBC approach of reconstructing rank-1 lattices, that
uses a single rank-1 lattice as sampling scheme,
one builds up finite sequences of rank-1 lattices which allow
for the exact reconstruction of trigonometric polynomials.
Numerical tests suggest significantly lower numbers $M$ of sampling nodes
that are required. In detail, numerical tests in \cite{Ka15} seem to
promise constant oversampling factors $M/|H_R^{d,0}|$.
Accordingly, the sampling rates could be possibly similar to those of sparse grids.

\section{Results for anisotropic mixed smoothness}
In this section we give an outlook on function spaces $\mathcal{H}^{\zb \alpha}_{\mathrm{mix}}(\T)$ where $\zb \alpha$
is a vector with first $\mu$ smallest smoothness directions, i.e.,
$$\frac{1}{2}<\alpha_1=\ldots=\alpha_{\mu}<\alpha_{\mu+1}\leq \ldots \leq \alpha_d.$$
\begin{definition}
	Let $\zb \alpha\in \R^d$ with positive entries. We define the Sobolev spaces with anisotropic mixed
smoothness $\zb \alpha$ as
	\begin{align*}
	\mathcal{H}^{\zb \alpha}_{\mathrm{mix}}(\T^d):=\left\{f\in L_2(\T^d)\colon\|f|\mathcal{H}^{\zb \alpha}_{\mathrm{mix}}(\T^d)\|^2:=\sum_{\zb k\in\Z^d}|\hat{f}_{\zb k}|^2\prod_{s=1}^d
	(1+|k_s|^2)^{\alpha_s}<\infty\right\}\,.
	\end{align*}
\end{definition}
Again, we want to study approximation by sampling along rank-1 lattices. Therefore we introduce new index sets,
so-called anisotropic hyperbolic crosses $H_R^{d,\zb \alpha}$ defined by
$$H_R^{d,\zb \alpha}:=\bigcup_{\zb j\in J_R^{d,\zb \alpha}}Q_{\zb j}$$
where
$$J_R^{d,\zb \alpha}:=\Big\{\zb j\in\N_0^d\colon \frac{1}{\alpha_1}\zb \alpha \cdot \zb j \leq R\Big\}.$$

\begin{Lemma}\label{lem:anisopcrosscard}
	Let $\zb \alpha\in\R^d$ with $0<\alpha_1=\ldots=\alpha_{\mu}<\alpha_{\mu+1}\leq\ldots\leq \alpha_{d}$. Then 
	$$|H^{d,\zb \alpha}_R|\asymp\sum_{\zb j\in J^{d,\zb \alpha}_R}2^{\|\zb j\|_1}\asymp 2^{R}R^{\mu-1}.$$
\end{Lemma}
\begin{proof}
	For the upper bound we refer to \cite[Chapt.\ 1., Lem.\ D]{Tem87}. For the lower bound we consider the subset
	$$J^{d,\zb \alpha}_{R,\mu}:=\{\zb j\in J^{d,\zb \alpha}_R: j_{\mu+1}=\ldots=j_d=0\}\subset J^{d,\zb \alpha}_{R}$$
	and obtain with the help of Lemma \ref{lem:card_HRdT}
	$$\sum_{\zb j\in J^{d,\zb \alpha}_R}2^{\|\zb j\|_1} \geq \sum_{\zb j\in J^{d,\zb \alpha}_{R,\mu}}2^{\|\zb j\|_1}\asymp\sum_{\zb j\in J^{\mu,0}_{R+c}}2^{\|\zb j\|_1}\gtrsim 2^{R}R^{\mu-1}.$$
\end{proof}

\begin{Lemma}\label{lem:anisobound_M}
	Let the refinement $R\ge1$, and the dimension $d\in\N$ with $d\ge2$, be given.
	Then there exists
	a reconstructing rank-1 lattice $\Lambda(\zb z,M)$ for $H_R^{d,{\zb \alpha}}$ which fulfills
	\begin{align*}
	2^RR^{\mu-1}\asymp |H^{d,\zb \alpha}_R| \leq M\lesssim
	2^{2R}R^{\mu-1}. 
	\end{align*}
\end{Lemma}
\begin{proof}
First, we show the embedding of the difference set $\mathcal{D}(H^{d,\zb \alpha}_R)\subset H^{d,\zb \alpha}_{2R+\|\zb \alpha\|_1}$.
Let $\zb k, \zb{k'}\in H_R^{d,\zb \alpha}$. Then there exist indices $\zb j,\zb{j'}\in J_R^{d,\zb \alpha}$ such that
$\zb k\in Q_{\zb j}$ and $\zb{k'}\in Q_{\zb{j'}}$. The difference $\zb k-\zb{k'}\in\mathcal{D}(H^{d,\zb \alpha}_R)$
and $\zb k-\zb{k'}\in Q_{\zb{\tilde j}}$ for an index $\zb{\tilde j}\in\N_0^d$. Next, we show $\zb \alpha\cdot\zb{\tilde j}\le2R+\|\zb \alpha\|_1$.
The differences $k_s-k'_s$ of one component of $\zb k$ and $\zb{k'}$ fulfill
$$
k_s-k'_s\in[-2^{j_s}-2^{j'_s},2^{j_s}+2^{j'_s}]\subset[-2^{\max(j_s,j'_s)+1},2^{\max(j_s,j'_s)+1}]=\bigcup_{t=0}^{\max(j_s,j'_s)+1}Q_t
$$
and we obtain ${\tilde j}_s\le\max(j_s,j'_s)+1\le j_s+j'_s+1$. This yields
$\zb\alpha\cdot\zb{\tilde j}\le\zb\alpha\cdot \zb{j}+\zb\alpha\cdot\zb{j'}+\|\zb \alpha\|_1\le 2R+\|\zb \alpha\|_1$
and consequently the embedding $\mathcal{D}(H^{d,\zb \alpha}_R)\subset H^{d,\zb \alpha}_{2R+\|\zb \alpha\|_1}$ holds.
Finally, the assertion is a consequence of Lemma \ref{lem:anisopcrosscard} and \cite[Corollary 3.4]{kaemmererdiss}.
\end{proof}

\begin{Remark}
The proof of Lemma \ref{lem:anisobound_M} referred here is based on an abstract result suitable for much more
general index sets than $H^{d,\zb \alpha}_R$. Similar to Lemma \ref{lem:universal_bounds_M} there should be also a
direct computation for counting the cardinality of the difference set $\mathcal{D}(H^{d,\zb \alpha}_R)$. We leave the
details to the interested reader. 
\end{Remark}

\begin{Lemma}\label{lem:anisosum}
	Let $\zb \alpha,\zb \gamma\in\R^d$ with
$\frac{1}{2}<\alpha_1=\gamma_1=\ldots=\alpha_{\mu}=\gamma_\mu<\alpha_{\mu+1}\leq\ldots\leq \alpha_{d}$ with
$\alpha_{\mu}<\gamma_s<\alpha_s$ for $s=\mu+1,\ldots,d$. Then it holds
	$$\sum_{\zb j\in \N_0^d\setminus J^{d,\zb \gamma}_R}2^{-(2\zb \alpha-\zb 1)\cdot \zb j}\lesssim 2^{-(2\alpha_1-1)R}R^{\mu-1}.$$
\end{Lemma}
\begin{proof}
	We start decomposing the sum. For technical reasons we introduce the notation
	$$P^{d, \zb \gamma}_R:= \Big\{\zb j\in \N_0^d: \frac{\gamma_s}{\gamma_1} j_s\leq R,s=1,\ldots,d\Big\}.$$
	Since $J^{d,\zb \gamma}_R \subset P^{d, \zb \gamma}_R$ we obtain
	\begin{align}
	\sum_{\zb j\in \N_0^d\setminus J^{d,\zb \gamma}_R}2^{-(2\zb \alpha-\zb 1)\cdot \zb j}= \sum_{\substack{\zb j\notin J^{d,\zb \gamma}_R\\\zb j\in P^{d, \zb \gamma}_R}}2^{-(2\zb \alpha-\zb 1)\cdot \zb j}+\sum_{\zb j\in \N_0^d\setminus P^{d, \zb \gamma}_R}2^{-(2\zb \alpha-\zb 1)\cdot \zb j}.\label{eq:decompsum}
	\end{align}
	We estimate the first summand in \eqref{eq:decompsum}
	\begin{align*}
	\sum_{\substack{\zb j\notin J^{d,\zb \gamma}_R\\\zb j\in P^{d, \zb \gamma}_R}}2^{-(2\zb \alpha- \zb 1)\cdot \zb j}=&\sum_{j_d=0}^{\frac{\gamma_1 R}{\gamma_d}}2^{-(2\alpha_d-1)j_d}\cdot\ldots\cdot\sum_{j_{\mu+1}=0}^{\frac{\gamma_1 R}{\gamma_{\mu+1}}}2^{-(2\alpha_{\mu+1}-1)j_{\mu+1}}\\&\cdot\sum_{j_{\mu}=0}^{\frac{\gamma_1R}{\gamma_{\mu}}}2^{-(2\alpha_{\mu}-1)j_{\mu}}\cdot\ldots\cdot \sum_{j_{2}=0}^{\frac{\gamma_1R}{\gamma_{2}}}2^{-(2\alpha_2-1)j_{2}} \sum_{j_{1}=\frac{\gamma_1 R-\sum_{s=2}^{d}\gamma_s j_s}{\gamma_1}}^{R}2^{-(2\alpha_1-1)j_{2}}\\
	&\lesssim \sum_{j_d=0}^{\frac{\gamma_1R}{\gamma_d}}2^{-(2\alpha_d-1)j_d}\cdot\ldots\cdot\sum_{j_{\mu+1}=0}^{\frac{\gamma_1R}{\gamma_{\mu+1}}}2^{-(2\alpha_{\mu+1}-1)j_{\mu+1}}\\&\cdot\sum_{j_{\mu}=0}^{\frac{\gamma_1R}{\gamma_{\mu}}}2^{-(2\alpha_{\mu}-1)j_{\mu}}\cdot\ldots\cdot \sum_{j_{2}=0}^{\frac{\gamma_1R}{\gamma_{2}}}2^{-(2\alpha_2-1)j_{2}} 2^{-(2\alpha_1-1)\frac{\gamma_1 R-\sum_{s=2}^{d}\gamma_s j_s}{\gamma_1}}.\\
	\end{align*}
Interchanging the order of multiplication yields
	\begin{align*}
	\sum_{\substack{\zb j\notin J^{d,\zb \gamma}_R\\\zb j\in P^{d, \zb \gamma}_R}}2^{-(2\zb \alpha- \zb 1)\cdot \zb j}&\lesssim 2^{-(2\alpha_1-1)R} \sum_{j_d=0}^{\infty}2^{-[(2\alpha_d-1)-(2\gamma_d-\frac{\gamma_d}{\alpha_1})]j_d}\cdot\ldots\\
	&\quad\cdot\sum_{j_{\mu+1}=0}^{\infty}2^{-[(2\alpha_{\mu+1}-1)-(2\gamma_{\mu+1}-\frac{\gamma_{\mu+1}}{\alpha_1})]j_{\mu+1}}
	\cdot\sum_{j_{\mu}=0}^{\frac{R}{\gamma_{\mu}-\varepsilon}}1\cdot\ldots\cdot \sum_{j_{2}=0}^{\frac{R}{\gamma_{2}-\varepsilon}}1\\
	&\lesssim 2^{-(2\alpha_1-1)R} R^{\mu-1}.
	\end{align*}
	The second summand in \eqref{eq:decompsum} can be trivially estimated by $\lesssim 2^{-(2\alpha_1-1)R}$.
\end{proof}
\begin{Theorem} \label{thm:vectorlinf}
	Let $\zb \alpha, \zb \gamma \in\R^d$ such that 
$$\frac{1}{2}<\alpha_1=\gamma_1=\ldots=\alpha_{\mu}=\gamma_\mu<\alpha_{\mu+1}\leq\ldots\leq \alpha_{d}$$ 
	and $$\alpha_1<\gamma_{s}<\alpha_s,\;s=\mu+1,\ldots,d,$$
	and the refinement $R\ge 1$, be given. In addition, we assume that $\Lambda(\zb z,M)$ is a reconstructing rank-1 lattice for $H_R^{d,\zb \gamma}$ 
	constructed by the CBC strategy \cite[Tab.\ 3.1]{kaemmererdiss}.
	We estimate the error of the sampling operator $\mathrm{Id}-S_{H_R^{d,\zb \gamma}}^{\Lambda(\zb z,M)}$ by
	\begin{align*}
	\|\mathrm{Id}-S_{H_R^{d,\zb \gamma}}^{\Lambda(\zb z,M)}|\mathcal{H}^{\zb \alpha}_{\mathrm{mix}}(\T^d)\rightarrow L_{\infty}(\T^d)\|
	&\lesssim 2^{-(\alpha_1-\frac{1}{2})R}
	R^{\frac{\mu-1}{2}}
	\\ &\lesssim M^{-(\alpha_1-\frac{1}{2})/2}
	(\log M)^{\frac{\mu-1}{2}(\alpha_1+\frac{1}{2})}.
	\end{align*}
\end{Theorem}
\begin{proof}
	We use the embedding $\mathcal{A}(\T^d)\hookrightarrow L_{\infty}(\T^d)$ and follow the estimation of Theorem  \ref{thm:alias_infty_err:dyadic} where we replace the weight $\prod_{s=1}^d(1+|k_s|^2)^{\alpha}$ by $\prod_{s=1}^d(1+|k_s|^2)^{\alpha_s}$. We obtain
	$$\|f-S_{H^{d,\zb \gamma}_R}^{\Lambda(\zb z,M)}f|L_{\infty}(\T^d)\|\lesssim \left(\sum_{\zb j\in \N_0^d\setminus J^{d,\zb \gamma}_R}2^{-(2\zb \alpha-\zb 1)\cdot \zb j}\right)^{\frac{1}{2}}\|f|\mathcal{H}^{\zb \alpha}_{\mathrm{mix}}(\T^d)\|.$$
	Applying Lemma \ref{lem:anisosum} yields
	$$ \|f-S_{H^{d,\zb \gamma}_R}^{\Lambda(\zb z,M)}f|L_{\infty}(\T^d)\|\lesssim 2^{-(\alpha_1-\frac{1}{2})R}R^{\frac{\mu-1}{2}}\|f|\mathcal{H}^{\zb \alpha}_{\mathrm{mix}}(\T^d)\|.$$
	Now the bound for the number of points in Lemma \ref{lem:anisobound_M} implies
	$$ \|f-S_{H^{d,\zb \gamma}_R}^{\Lambda(\zb z,M)}f|L_{\infty}(\T^d)\|\lesssim M^{-(\alpha_1-\frac{1}{2})/2}	(\log M)^{\frac{\mu-1}{2}(\alpha_1+\frac{1}{2})}\|f|\mathcal{H}^{\zb \alpha}_{\mathrm{mix}}(\T^d)\| .$$
	That proves the claim.
\end{proof}
\begin{Remark}
	Comparing the last result with the results obtained in Proposition \ref{thm:alias_err_infty:dyadic} we
recognize that there is only the exponent $\mu-1$ instead of $d-1$ in the logarithm of the
error term with $\mu<d$. Especially in the case $\mu=1$ the logarithm completely vanishes.
Similar effects were also observed for sparse grids and general linear approximation, cf. \cite{DD15,Tem87}.
\end{Remark}

\section{Numerical results}

In this section, we numerically investigate the sampling rates for different types of rank-1 lattices $\Lambda(\zb
z,M)$ when sampling the scaled periodized (tensor product) kink function
\begin{equation} \label{equ:kink}
g(\boldx):= \prod_{t=1}^d \left(\frac{5^{3/4} 15}{4 \sqrt{3}} \max\left\{ \frac{1}{5} - \left(x_t-\frac{1}{2}\right)^2,0 \right\} \right), \quad \boldx:=(x_1,\ldots,x_d)^\top\in\T^d,
\end{equation}
similar to \cite{HiMaOeUl15}.
We remark that $g\in\mathcal{H}^{3/2-\varepsilon}_\mathrm{mix}(\T^d)$, $\varepsilon>0$, and $\|g|L_2(\T^d)\|=1$.

For the fast approximate reconstruction, Algorithm~\ref{alg:r1l_reconstruction_fct} can be used. This algorithm applies a single one-dimensional fast Fourier transform (FFT) on the function samples and performs a simple index transform. As input parameter a reconstructing rank-1 lattice $\Lambda(\boldz,M)$ is required, which may be easily searched for by means of the CBC strategy \cite[Tab.\ 3.1]{kaemmererdiss}.

\begin{algorithm}[ht]
\caption{Fast approximate reconstruction of a function 
$f\in\mathcal{H}^{\alpha,\beta}(\T^d)$
from sampling values on a reconstructing rank-1 lattice $\Lambda(\boldz,M)$ using a single one-dimensional FFT, see \cite[Algorithm~1]{KaPoVo13}.}\label{alg:r1l_reconstruction_fct}
  \begin{tabular}{p{2.1cm}p{4.5cm}p{7.7cm}}
    Input:      
                & $ I\subset\Z^d$ \hfill & frequency index set of finite cardinality  \\
                & $\Lambda(\boldz,M)$ & reconstructing rank-1 lattice for $I$ of size $M$ with generating vector $\boldz\in\Z^d$\\
                & $\zb{f}=\left(f\left(\frac{j\zb z}{M}\bmod{\zb 1}\right)\right)_{j=0}^{M-1}$ & samples of $f\in\mathcal{H}^{\alpha,\beta}(\T^d)$ on $\Lambda(\boldz,M)$
  \end{tabular}
  \begin{algorithmic}
        \STATE $\zb {\hat{a}}:=\mathrm{FFT\_1D}(\zb f)$
        \FOR{\textbf{each} $\boldk\in I$}
        \STATE  $\hat{f}_{\zb k}^{\Lambda(\zb z,M)}:=\frac{1}{M} \, \hat{a}_{\boldk\cdot\boldz\bmod{M}}$
        \ENDFOR
  \end{algorithmic}
  \begin{tabular}{p{2.1cm}p{4.5cm}p{7.7cm}}
    Output: & $\hat{f}_{\zb k}^{\Lambda(\zb z,M)}$  
    & Fourier coefficients of the approximation $S_I^{\Lambda(\zb z,M)}f$ as defined in~\eqref{eqn:sampling_operator} \\
    Complexity: &  \multicolumn{2}{l}{$\mathcal{O}\left(M\log M + d\vert I\vert\right)$}
  \end{tabular}
\end{algorithm}

\subsection{Hyperbolic cross index sets}
\label{sec:numberics:kink:hc}

First, we build reconstructing rank-1 lattices for the hyperbolic cross index sets $H^{d,0}_R$ in the cases $d=2,3,4$ with various refinements $R\in\N_0$
using the CBC strategy \cite[Tab.\ 3.1]{kaemmererdiss}. Then, we apply the sampling operators $S^{\Lambda(\zb
z,M)}_{H^{d,0}_R}$ on the kink function $g$ using Algorithm~\ref{alg:r1l_reconstruction_fct}.
The resulting sampling errors $\|g-S^{\Lambda(\zb z,M)}_{H^{d,0}_R}g|L_2(\T^d)\|$ are shown in Figure~\ref{figure:numerics:kink:hc} and \ref{figure:numerics:kink:hc2} denoted by ``CBC hc''.
The corresponding theoretical upper bounds for the sampling rates from Table~\ref{tab:mainres:overview_general_r1l}, which are (almost)
$M^{-\frac{1}{2}\cdot\frac{3}{2}}(\log M)^{\frac{d-2}{2}\cdot\frac{3}{2}+\frac{d-1}{2}}$, are also depicted.
Additionally in the two-dimensional case, we consider the Fibonacci lattices from Section~\ref{sec:2d} as well as special Korobov lattices $$\Lambda((1,\lceil 3\cdot 2^{R-2}\rceil)^\top, \lceil(1+3\cdot 2^{R-2})\cdot 2^{R-1}\rceil )$$ from \cite{KaKuPo10}.
The corresponding sampling errors are denoted by ``Fibonacci hc'' and ``Korobov hc'' in Figure~\ref{figure:numerics:kink:hc}.
We observe that in all considered cases, the sampling errors decay approximately as fast as the theoretical upper bound implies.
In Figure~\ref{figure:numerics:kink:cbc_hc_scaled}, we investigate the logarithmic factors in more detail.
Assuming that the sampling error $\|g-S^{\Lambda(\zb z,M)}_{H^{d,0}_R}g|L_2(\T^d)\|$ nearly decays like
$M^{-\frac{1}{2}\cdot\frac{3}{2}}(\log M)^{\frac{d-2}{2}\cdot\frac{3}{2}+\frac{d-1}{2}}$,
we consider its scaled version $$\|g-S^{\Lambda(\zb z,M)}_{H^{d,0}_R}g|L_2(\T^d)\| /[ M^{-\frac{1}{2}\cdot\frac{3}{2}}(\log M)^{\frac{d-2}{2}\cdot\frac{3}{2}+\frac{d-1}{2}}].$$
Obviously, if the scaled error decays exactly like the given rate, then the plot should be (approximately) a horizontal line.
In the plot in Figure~\ref{figure:numerics:kink:hc_scaled_d2:scaledlog} for the two-dimensional case, this is almost the case for all three types of lattices. The scaled errors $\|g-S^{\Lambda(\zb z,M)}_{H^{d,0}_R}g|L_2(\T^d)\| \cdot M^{1.5/2} \cdot (\log M)^{-1/2}$ seem to decay slightly but the errors in Figure \ref{figure:numerics:kink:hc_scaled_d2:scaledmain}, that are scaled without the logarithmic factor, grow slightly. We interpret this observation as an indication that there is some logarithmic dependence in the error rate.
Moreover, for the reconstructing rank-1 lattices built using the CBC strategy \cite[Tab.\ 3.1]{kaemmererdiss},
the scaled errors in the cases $d=3$, $d=4$, and $d=5$ behave similarly as in the two-dimensional case, see Figure~\ref{figure:numerics:kink:hc_scaled_d234}.

\begin{figure}[ht]
\centering
\subfloat[$d=2$]{\label{figure:numerics:kink:hc:d2}
\begin{tikzpicture}
  \begin{loglogaxis}[enlargelimits=false,xmin=1e1,xmax=1e9,ymin=1e-7,ymax=2e0,ytick={1e-6,1e-4,1e-2,1},height=0.34\textwidth, width=0.92\textwidth, grid=major, xlabel={Number of sampling points $M$}, ylabel={$\|g-S^{\Lambda(\zb z,M)}_{H^{d,0}_R}g|L_2(\T^d)\|$},
  legend style={at={(0.5,1.05)}, anchor=south,legend columns=-1, font=\footnotesize},
  ]
\addplot[blue,mark=*,mark size=2] coordinates { 
 (2,1.000e+00) (5,5.424e-01) (13,2.825e-01) (47,1.431e-01) (163,4.521e-02) (593,1.259e-02) (2243,5.884e-03) (8641,2.152e-03) (33797,7.110e-04) (133379,2.698e-04) (529411,9.918e-05) (2108429,3.534e-05) (8413187,1.286e-05) (33607697,4.669e-06) (134332421,1.691e-06) 
};
\addlegendentry{CBC hc}
\addplot[red,mark=square*,mark size=2] coordinates { 
(1,1.000e+00) (2,7.322e-01) (3,6.656e-01) (5,5.424e-01) (8,5.190e-01) (13,2.497e-01) (21,2.460e-01) (34,1.413e-01) (55,1.410e-01) (89,1.407e-01) (144,4.486e-02) (233,4.443e-02) (377,4.440e-02) (610,1.201e-02) (987,1.195e-02) (1597,5.718e-03) (2584,5.713e-03) (4181,5.710e-03) (6765,5.709e-03) (10946,2.110e-03) (17711,2.110e-03) (28657,6.923e-04) (46368,6.918e-04) (75025,6.918e-04) (121393,6.918e-04) (196418,2.645e-04) (317811,2.644e-04) (514229,9.692e-05) (832040,9.692e-05) (1346269,9.692e-05) (2178309,9.692e-05) (3524578,3.463e-05) (5702887,3.463e-05) (9227465,1.264e-05) (14930352,1.264e-05) (24157817,4.588e-06) (39088169,4.588e-06) (63245986,4.588e-06) (102334155,4.588e-06) (165580141,1.659e-06) (267914296,1.659e-06) (433494437,6.009e-07) 
 };
\addlegendentry{Fibonacci hc}
\addplot[brown,mark=diamond*,mark size=2.5] coordinates { 
 (1,1.000e+00) (3,6.656e-01) (8,3.022e-01) (28,1.452e-01) (104,4.602e-02) (400,1.285e-02) (1568,5.895e-03) (6208,2.157e-03) (24704,7.106e-04) (98560,2.701e-04) (393728,9.928e-05) (1573888,3.535e-05) (6293504,1.287e-05) (25169920,4.669e-06) (100671488,1.690e-06) (402669568,6.102e-07) 
 };
\addlegendentry{Korobov hc}
\addplot[black,domain=1:1e9,samples=100,dashed] {2.5*x^(-1.5/2)*ln(x)^((2-2)*1.5/2+(2-1)/2)};
\addlegendentry{$\text{2.5}\cdot M^{-\frac{3}{4}} (\log M)^{\frac{1}{2}}$}
  \end{loglogaxis}
\end{tikzpicture}
}
\hfill
\subfloat[$d=3$]{
\begin{tikzpicture}
  \begin{loglogaxis}[enlargelimits=false,xmin=1e1,xmax=1e9,ymin=3e-5,ymax=2e0,ytick={1e-4,1e-3,1e-2,1e-1,1},height=0.34\textwidth, width=0.92\textwidth, grid=major, xlabel={Number of sampling points $M$}, ylabel={$\|g-S^{\Lambda(\zb z,M)}_{H^{d,0}_R}g|L_2(\T^d)\|$},
  legend style={at={(0.5,1.05)}, anchor=south,legend columns=-1, font=\footnotesize},
  ]
\addplot[blue,mark=*,mark size=2] coordinates { 
 (2,1.162e+00) (5,7.310e-01) (29,4.275e-01) (149,2.543e-01) (787,1.162e-01) (3877,5.798e-02) (18661,2.145e-02) (87403,7.024e-03) (400523,2.980e-03) (1805053,1.168e-03) (8029361,4.127e-04) (35344873,1.569e-04) (154231373,5.955e-05)
 };
\addlegendentry{CBC hc}
\addplot[black,domain=1e1:1e9,samples=100,dashed] {1.5*x^(-1.5/2)*ln(x)^((3-2)*1.5/2+(3-1)/2)};
\addlegendentry{$\text{1.5}\cdot M^{-\frac{3}{4}}(\log M)^{\frac{7}{4}}$}
  \end{loglogaxis}
\end{tikzpicture}
}
\hfill
\subfloat[$d=4$]{
\begin{tikzpicture}
  \begin{loglogaxis}[enlargelimits=false,xmin=1e1,xmax=1e9,ymin=4e-4,ymax=2e0,ytick={1e-3,1e-2,1e-1,1},height=0.34\textwidth, width=0.92\textwidth, grid=major, xlabel={Number of sampling points $M$}, ylabel={$\|g-S^{\Lambda(\zb z,M)}_{H^{d,0}_R}g|L_2(\T^d)\|$},
  legend style={at={(0.5,1.05)}, anchor=south,legend columns=-1, font=\footnotesize},
  ]
\addplot[blue,mark=*,mark size=2] coordinates {
(2,1.840e+00) (7,7.700e-01) (53,5.093e-01) (389,3.535e-01) (2579,1.951e-01) (15401,1.123e-01) (87337,5.374e-02) (472907,2.502e-02) (2473697,1.006e-02) (12560731,3.698e-03) (62357357,1.515e-03) (303057067,6.083e-04)
};
\addlegendentry{CBC hc}
\addplot[black,domain=1e1:1e9,samples=100,dashed] {0.5*x^(-1.5/2)*ln(x)^((4-2)*1.5/2+(4-1)/2)};
\addlegendentry{$\text{5e-1}\cdot M^{-\frac{3}{4}}(\log M)^{3}$}
  \end{loglogaxis}
\end{tikzpicture}
}
\caption{$L_2(\T^d)$ sampling error and number of sampling points for the approximation of the kink function $g$ from \eqref{equ:kink}.}
\label{figure:numerics:kink:hc}
\end{figure}
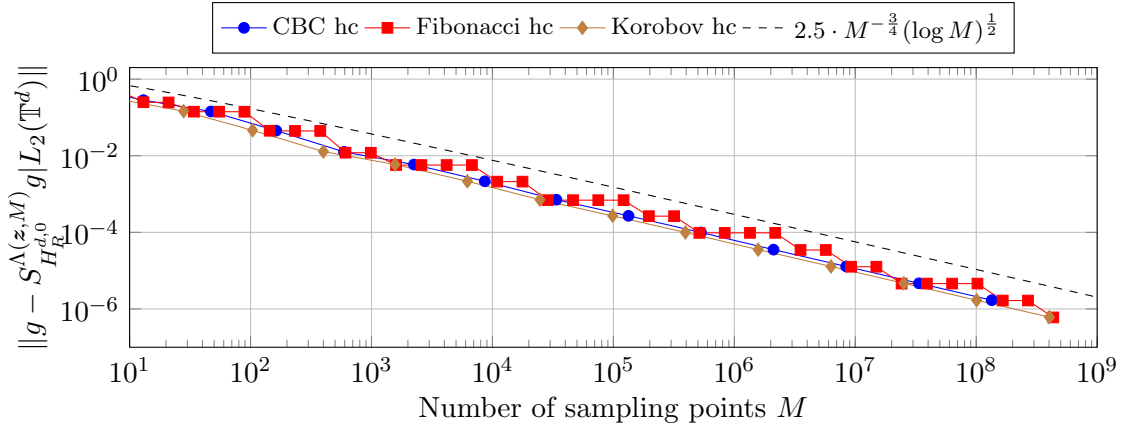
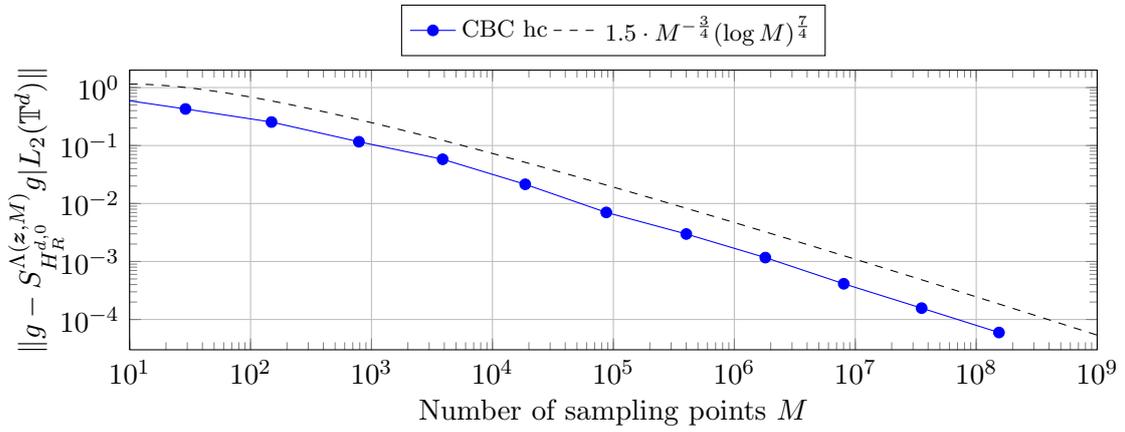
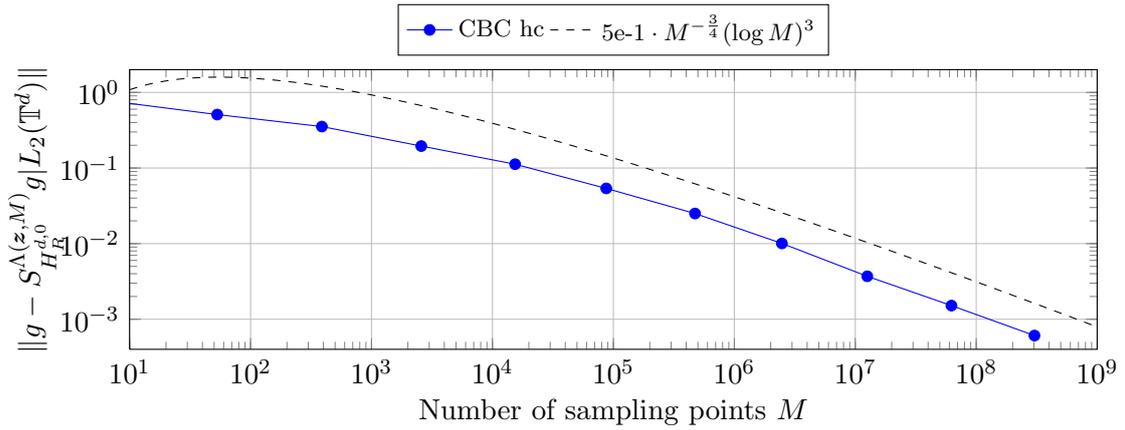

\begin{figure}[ht]
\centering
\subfloat[$d=5$]{
\begin{tikzpicture}
  \begin{loglogaxis}[enlargelimits=false,xmin=1e1,xmax=1e9,ymin=3e-3,ymax=3e0,ytick={1e-3,1e-2,1e-1,1},height=0.34\textwidth, width=0.92\textwidth, grid=major, xlabel={Number of sampling points $M$}, ylabel={$\|g-S^{\Lambda(\zb z,M)}_{H^{d,0}_R}g|L_2(\T^d)\|$},
  legend style={at={(0.5,1.05)}, anchor=south,legend columns=-1, font=\footnotesize},
  ]
\addplot[blue,mark=*,mark size=2] coordinates {
(2,2.841e+00) (7,8.391e-01) (89,6.155e-01) (863,4.499e-01) (6863,2.766e-01) (48073,1.732e-01) (311279,9.619e-02) (1893643,5.190e-02) (11018881,2.483e-02) (61634213,1.126e-02) (335314513,4.717e-03)
};
\addlegendentry{CBC hc}
\addplot[black,domain=1e1:1e9,samples=100,dashed] {0.1*x^(-1.5/2)*ln(x)^((5-2)*1.5/2+(5-1)/2)};
\addlegendentry{$\text{1e-1}\cdot M^{-\frac{3}{4}}(\log M)^{\frac{17}{4}}$}
  \end{loglogaxis}
\end{tikzpicture}
}
\hfill
\subfloat[$d=6$]{
\begin{tikzpicture}
  \begin{loglogaxis}[enlargelimits=false,xmin=1e1,xmax=1e9,ymin=1e-2,ymax=4e0,ytick={1e-2,1e-1,1},height=0.34\textwidth, width=0.92\textwidth, grid=major, xlabel={Number of sampling points $M$}, ylabel={$\|g-S^{\Lambda(\zb z,M)}_{H^{d,0}_R}g|L_2(\T^d)\|$},
  legend style={at={(0.5,1.05)}, anchor=south,legend columns=-1, font=\footnotesize},
  ]
\addplot[blue,mark=*,mark size=2] coordinates {
(2,4.290e+00) (11,8.843e-01) (149,7.081e-01) (1709,5.252e-01) (16087,3.537e-01) (128969,2.377e-01) (937537,1.458e-01) (6312689,8.639e-02) (40272361,4.671e-02) (244936711,2.415e-02)
};
\addlegendentry{CBC hc}
\addplot[black,domain=1e1:1e9,samples=100,dashed] {0.01*x^(-1.5/2)*ln(x)^((6-2)*1.5/2+(6-1)/2)};
\addlegendentry{$\text{1e-2}\cdot M^{-\frac{3}{4}}(\log M)^{\frac{11}{2}}$}
  \end{loglogaxis}
\end{tikzpicture}
}
\hfill
\subfloat[$d=7$]{
\begin{tikzpicture}
  \begin{loglogaxis}[enlargelimits=false,xmin=1e1,xmax=1e9,ymin=3e-2,ymax=3e0,ytick={1e-1,1},height=0.34\textwidth, width=0.92\textwidth, grid=major, xlabel={Number of sampling points $M$}, ylabel={$\|g-S^{\Lambda(\zb z,M)}_{H^{d,0}_R}g|L_2(\T^d)\|$},
  legend style={at={(0.5,1.05)}, anchor=south,legend columns=-1, font=\footnotesize},
  ]
\addplot[blue,mark=*,mark size=2] coordinates {
(2,6.380e+00) (11,9.099e-01) (223,7.655e-01) (3137,6.011e-01) (34283,4.286e-01) (310861,3.059e-01) (2510843,2.000e-01) (18545321,1.271e-01) (128547037,7.504e-02) (843436063,4.2478e-02)
};
\addlegendentry{CBC hc}
\addplot[black,domain=1e1:1e9,samples=100,dashed] {0.0006*x^(-1.5/2)*ln(x)^((7-2)*1.5/2+(7-1)/2)};
\addlegendentry{$\text{6e-4}\cdot M^{-\frac{3}{4}}(\log M)^{\frac{27}{4}}$}
  \end{loglogaxis}
\end{tikzpicture}
}
\caption{$L_2(\T^d)$ sampling error and number of sampling points for the approximation of the kink function $g$ from \eqref{equ:kink}.}
\label{figure:numerics:kink:hc2}
\end{figure}

\begin{figure}[ht]
\centering
\subfloat[$d=2$]{\label{figure:numerics:kink:hc_scaled_d2:scaledlog}
\begin{tikzpicture}
\begin{loglogaxis}[enlargelimits=false,height=0.38\textwidth, width=0.95\textwidth, grid=major, xlabel={Number of sampling points $M$}, ylabel={$\mathrm{err} \cdot M^{1.5/2} \cdot (\log M)^{-1/2}$},xmin=1e1,xmax=6e8,ymin=1e-1,ymax=1e1,
  legend style={at={(0.5,1.05)}, anchor=south,legend columns=-1, font=\footnotesize},
  ]
\addplot[blue,mark=*,mark size=2] table[row sep=\\,x=M, y expr = { \thisrow{error} * (\thisrow{M})^(3/4) * ln(\thisrow{M})^(-1/2) }] {
M error\\
13 2.825e-01\\
47 1.431e-01\\
163 4.521e-02\\
593 1.259e-02\\
2243 5.884e-03\\
8641 2.152e-03\\
33797 7.110e-04\\
133379 2.698e-04\\
529411 9.918e-05\\
2108429 3.534e-05\\
8413187 1.286e-05\\
33607697 4.669e-06\\
134332421 1.691e-06\\
};
\addlegendentry{CBC hc}
\addplot[red,mark=square*,mark size=2] table[row sep=\\,x=M, y expr = { \thisrow{error} * (\thisrow{M})^(3/4) * ln(\thisrow{M})^(-1/2) }] {
M error\\
13 2.497e-01\\
21 2.460e-01\\
34 1.413e-01\\
55 1.410e-01\\
89 1.407e-01\\
144 4.486e-02\\
233 4.443e-02\\
377 4.440e-02\\
610 1.201e-02\\
987 1.195e-02\\
1597 5.718e-03\\
2584 5.713e-03\\
4181 5.710e-03\\
6765 5.709e-03\\
10946 2.110e-03\\
17711 2.110e-03\\
28657 6.923e-04\\
46368 6.918e-04\\
75025 6.918e-04\\
121393 6.918e-04\\
196418 2.645e-04\\
317811 2.644e-04\\
514229 9.692e-05\\
832040 9.692e-05\\
1346269 9.692e-05\\
2178309 9.692e-05\\
3524578 3.463e-05\\
5702887 3.463e-05\\
9227465 1.264e-05\\
14930352 1.264e-05\\
24157817 4.588e-06\\
39088169 4.588e-06\\
63245986 4.588e-06\\
102334155 4.588e-06\\
165580141 1.659e-06\\
267914296 1.659e-06\\
433494437 6.009e-07\\
};
\addlegendentry{Fibonacci hc}
\addplot[brown,mark=diamond*,mark size=2.5] table[row sep=\\,x=M, y expr = { \thisrow{error} * (\thisrow{M})^(3/4) * ln(\thisrow{M})^(-1/2) }] {
M error\\
28 1.452e-01\\
104 4.602e-02\\
400 1.285e-02\\
1568 5.895e-03\\
6208 2.157e-03\\
24704 7.106e-04\\
98560 2.701e-04\\
393728 9.928e-05\\
1573888 3.535e-05\\
6293504 1.287e-05\\
25169920 4.669e-06\\
100671488 1.690e-06\\
402669568 6.102e-07\\
};
\addlegendentry{Korobov hc}
\end{loglogaxis}
\end{tikzpicture}
}
\hfill
\subfloat[$d=2$]{\label{figure:numerics:kink:hc_scaled_d2:scaledmain}
\begin{tikzpicture}
\begin{loglogaxis}[enlargelimits=false,height=0.38\textwidth, width=0.95\textwidth, grid=major, xlabel={Number of sampling points $M$}, ylabel={$\mathrm{err} \cdot M^{1.5/2} $},xmin=1e1,xmax=6e8,ymin=1e-1,ymax=1e1,
  legend style={at={(0.5,1.05)}, anchor=south,legend columns=-1, font=\footnotesize},
  ]
\addplot[blue,mark=*,mark size=2] table[row sep=\\,x=M, y expr = { \thisrow{error} * (\thisrow{M})^(3/4)  }] {
M error\\
13 2.825e-01\\
47 1.431e-01\\
163 4.521e-02\\
593 1.259e-02\\
2243 5.884e-03\\
8641 2.152e-03\\
33797 7.110e-04\\
133379 2.698e-04\\
529411 9.918e-05\\
2108429 3.534e-05\\
8413187 1.286e-05\\
33607697 4.669e-06\\
134332421 1.691e-06\\
};
\addlegendentry{CBC hc}
\addplot[red,mark=square*,mark size=2] table[row sep=\\,x=M, y expr = { \thisrow{error} * (\thisrow{M})^(3/4) }] {
M error\\
13 2.497e-01\\
21 2.460e-01\\
34 1.413e-01\\
55 1.410e-01\\
89 1.407e-01\\
144 4.486e-02\\
233 4.443e-02\\
377 4.440e-02\\
610 1.201e-02\\
987 1.195e-02\\
1597 5.718e-03\\
2584 5.713e-03\\
4181 5.710e-03\\
6765 5.709e-03\\
10946 2.110e-03\\
17711 2.110e-03\\
28657 6.923e-04\\
46368 6.918e-04\\
75025 6.918e-04\\
121393 6.918e-04\\
196418 2.645e-04\\
317811 2.644e-04\\
514229 9.692e-05\\
832040 9.692e-05\\
1346269 9.692e-05\\
2178309 9.692e-05\\
3524578 3.463e-05\\
5702887 3.463e-05\\
9227465 1.264e-05\\
14930352 1.264e-05\\
24157817 4.588e-06\\
39088169 4.588e-06\\
63245986 4.588e-06\\
102334155 4.588e-06\\
165580141 1.659e-06\\
267914296 1.659e-06\\
433494437 6.009e-07\\
};
\addlegendentry{Fibonacci hc}
\addplot[brown,mark=diamond*,mark size=2.5] table[row sep=\\,x=M, y expr = { \thisrow{error} * (\thisrow{M})^(3/4) }] {
M error\\
28 1.452e-01\\
104 4.602e-02\\
400 1.285e-02\\
1568 5.895e-03\\
6208 2.157e-03\\
24704 7.106e-04\\
98560 2.701e-04\\
393728 9.928e-05\\
1573888 3.535e-05\\
6293504 1.287e-05\\
25169920 4.669e-06\\
100671488 1.690e-06\\
402669568 6.102e-07\\
};
\addlegendentry{Korobov hc}
\end{loglogaxis}
\end{tikzpicture}
}
\hfill
\subfloat[$d=2,3,4,5$]{\label{figure:numerics:kink:hc_scaled_d234}
\begin{tikzpicture}
\begin{loglogaxis}[enlargelimits=false,height=0.38\textwidth, width=0.95\textwidth, grid=major, xlabel={Number of sampling points $M$}, ylabel={$\mathrm{err} \cdot M^{1.5/2} \cdot (\log M)^{-\frac{d-2}{2}\cdot\frac{3}{2}-\frac{d-1}{2}}$},xmin=1e1,xmax=6e8,ymin=1e-2,ymax=1e1,
  legend style={at={(0.5,1.05)}, anchor=south,legend columns=-1, font=\footnotesize},
  ]
\addplot[blue,mark=*,mark size=2] table[row sep=\\,x=M, y expr = { \thisrow{error} * (\thisrow{M})^(3/4) * ln(\thisrow{M})^(-1/2) }] {
M error\\
13 2.825e-01\\
47 1.431e-01\\
163 4.521e-02\\
593 1.259e-02\\
2243 5.884e-03\\
8641 2.152e-03\\
33797 7.110e-04\\
133379 2.698e-04\\
529411 9.918e-05\\
2108429 3.534e-05\\
8413187 1.286e-05\\
33607697 4.669e-06\\
134332421 1.691e-06\\
};
\addlegendentry{CBC hc $d=2$}
\addplot[red,mark=square*,mark size=2] table[row sep=\\,x=M, y expr = { \thisrow{error} * (\thisrow{M})^(3/4) * ln(\thisrow{M})^(-(3-2)/2 * 3/2 - (3-1)/2) }] {
M error\\
29 4.275e-01\\
149 2.543e-01\\
787 1.162e-01\\
3877 5.798e-02\\
18661 2.145e-02\\
87403 7.024e-03\\
400523 2.980e-03\\
1805053 1.168e-03\\
8029361 4.127e-04\\
35344873 1.569e-04\\
154231373 5.955e-05\\
};
\addlegendentry{CBC hc $d=3$}
\addplot[brown,mark=diamond*,mark size=2.5] table[row sep=\\,x=M, y expr = { \thisrow{error} * (\thisrow{M})^(3/4) * ln(\thisrow{M})^(-(4-2)/2 * 3/2 - (4-1)/2) }] {
M error\\
53 5.093e-01\\
389 3.535e-01\\
2579 1.951e-01\\
15401 1.123e-01\\
87337 5.374e-02\\
472907 2.502e-02\\
2473697 1.006e-02\\
12560731 3.698e-03\\
62357357 1.515e-03\\
303057067 6.083e-04\\
};
\addlegendentry{CBC hc $d=4$}
\addplot[black,mark=x,mark size=2.5] table[row sep=\\,x=M, y expr = { \thisrow{error} * (\thisrow{M})^(3/4) * ln(\thisrow{M})^(-(5-2)/2 * 3/2 - (5-1)/2) }] {
M error\\
89 6.155e-01\\
863 4.499e-01\\
6863 2.766e-01\\
48073 1.732e-01\\
311279 9.619e-02\\
1893643 5.190e-02\\
11018881 2.483e-02\\
61634213 1.126e-02\\
335314513 4.717e-03\\
};
\addlegendentry{CBC hc $d=5$}
\end{loglogaxis}
\end{tikzpicture}
}
\caption{Scaled $L_2(\T^d)$ sampling error and number of sampling points for the approximation of the kink function $g$ from \eqref{equ:kink},
where $\mathrm{err}:=\|g-S^{\Lambda(\zb z,M)}_{H^{d,0}_R}g|L_2(\T^d)\|$.}
\label{figure:numerics:kink:cbc_hc_scaled}
\end{figure}

\FloatBarrier

\subsection{$\ell_\infty$-ball index sets}
\begin{sloppypar}
Next, we use the lattices from Section~\ref{sec:numberics:kink:hc} in the two-dimensional case, but instead of hyperbolic cross index sets $H^{2,0}_R$,
we are going to use the $\ell_\infty$-ball index sets
$I^2_N:=\left\{ -\left\lceil\frac{N-2}{2}\right\rceil,\ldots, \left\lceil\frac{N-1}{2}\right\rceil \right\}^2$, $N\in\N$.
For each of the rank-1 lattices $\Lambda(\zb z,M)$ generated in Section~\ref{sec:numberics:kink:hc},
we determine the largest refinement $N\in\N$ such that the reconstruction property \eqref{eqn:reco_prop1} is still fulfilled for the $\ell_\infty$-ball $I^2_N$.
Then, we apply each sampling operator $S^{\Lambda(\zb z,M)}_{I^2_N}$ on the kink function $g$ from~\eqref{equ:kink}.
The resulting sampling errors are depicted in Figure~\ref{figure:numerics:kink:linf_d2}, where the
errors for the CBC, Fibonacci and Korobov rank-1 lattices are denoted by
``CBC $\ell_\infty$-ball'', ``Fibonacci $\ell_\infty$-ball'' and ``Korobov $\ell_\infty$-ball'', respectively.
We observe that the $L_2(\T^d)$ sampling errors decay approximately as the rate $M^{-\frac{3}{4}}$ as expected.
In more detail, this behaviour may be seen in the scaled error plot in Figure~\ref{figure:kink:linf_scaled}.
\end{sloppypar}

\begin{figure}[ht]
\centering
\begin{tikzpicture}
  \begin{loglogaxis}[enlargelimits=false,xmin=0.5,xmax=1e9,ymin=1e-7,ymax=2e0,ytick={1e-6,1e-4,1e-2,1},height=0.38\textwidth, width=0.94\textwidth, grid=major, xlabel={Number of sampling points $M$}, ylabel={$\|g-S^{\Lambda(\zb z,M)}_{I^2_N}g|L_2(\T^d)\|$}, xmin=1e1,
  legend style={at={(0.5,1.05)}, anchor=south,legend columns=-1, font=\footnotesize},
  ]
\addplot[blue,mark=*,mark size=2] coordinates { 
 (2,1.000e+00) (5,6.287e-01) (13,1.835e-01) (47,5.198e-02) (163,2.521e-02) (593,7.517e-03) (2243,3.274e-03) (8641,1.243e-03) (33797,4.215e-04) (133379,1.546e-04) (529411,5.635e-05) (2108429,1.948e-05) (8413187,6.947e-06) (33607697,2.498e-06) (134332421,8.864e-07) 
};
\addlegendentry{CBC $\ell_\infty$-ball}
\addplot[red,mark=square*,mark size=2] coordinates { 
(1,1.000e+00) (2,7.322e-01) (3,7.274e-01) (5,6.287e-01) (8,5.158e-01) (13,1.685e-01) (21,1.585e-01) (34,5.964e-02) (55,5.314e-02) (89,2.388e-02) (144,2.158e-02) (233,1.847e-02) (377,1.790e-02) (610,7.939e-03) (987,7.659e-03) (1597,4.402e-03) (2584,4.273e-03) (4181,2.040e-03) (6765,1.974e-03) (10946,9.793e-04) (17711,9.469e-04) (28657,4.636e-04) (46368,4.474e-04) (75025,2.182e-04) (121393,2.100e-04) (196418,1.087e-04) (317811,1.048e-04) (514229,5.211e-05) (832040,5.021e-05) (1346269,2.541e-05) (2178309,2.449e-05) (3524578,1.235e-05) (5702887,1.190e-05) (9227465,6.009e-06) (14930352,5.792e-06) (24157817,2.925e-06) (39088169,2.820e-06) (63245986,1.419e-06) (102334155,1.368e-06) (165580141,6.900e-07) (267914296,6.652e-07) (433494437,3.352e-07) 
};
\addlegendentry{Fibonacci $\ell_\infty$-ball}
\addplot[brown,mark=diamond*,mark size=2.5] coordinates { 
(1,1.000e+00) (3,7.274e-01) (8,5.148e-01) (28,1.181e-01) (104,2.457e-02) (400,1.528e-02) (1568,5.131e-03) (6208,1.584e-03) (24704,5.992e-04) (98560,2.079e-04) (393728,7.520e-05) (1573888,2.602e-05) (6293504,9.097e-06) (25169920,3.286e-06) (100671488,1.149e-06) (402669568,4.070e-07) 
};
\addlegendentry{Korobov $\ell_\infty$-ball}
\addplot[black,domain=1:1e9,samples=100] {2.5*x^(-1.5/2)};
\addlegendentry{$\text{2.5}\cdot M^{-\frac{3}{4}} $}
\addplot[black,domain=1:1e9,samples=100,dashed] {1.7*x^(-1.5/2)*ln(x)^((2-2)*1.5/2+(2-1)/2)};
\addlegendentry{$\text{1.7}\cdot M^{-\frac{3}{4}} (\log M)^{\frac{1}{2}}$}
  \end{loglogaxis}
\end{tikzpicture}
\caption{$L_2(\T^2)$ sampling error and number of sampling points for the approximation of the kink function $g$ from \eqref{equ:kink}.}
\label{figure:numerics:kink:linf_d2}

\begin{tikzpicture}
\begin{loglogaxis}[enlargelimits=false,height=0.38\textwidth, width=0.95\textwidth, grid=major, xlabel={Number of sampling points $M$}, ylabel={$\mathrm{err} \cdot M^{1.5/2}$}, xmin=1e1,ymin=1e-1,ymax=1e1,
  legend style={at={(0.5,1.05)}, anchor=south,legend columns=-1, font=\footnotesize},
  ]
\addplot[blue,mark=*,mark size=2] table[row sep=\\,x=M, y expr = { \thisrow{error} * (\thisrow{M})^(3/4) }] {
M error\\
2 1.000e+00\\
5 6.287e-01\\
13 1.835e-01\\
47 5.198e-02\\
163 2.521e-02\\
593 7.517e-03\\
2243 3.274e-03\\
8641 1.243e-03\\
33797 4.215e-04\\
133379 1.546e-04\\
529411 5.635e-05\\
2108429 1.948e-05\\
8413187 6.947e-06\\
33607697 2.498e-06\\
134332421 8.864e-07\\
};
\addlegendentry{CBC $\ell_\infty$-ball}
\addplot[red,mark=square*,mark size=2] table[row sep=\\,x=M, y expr = { \thisrow{error} * (\thisrow{M})^(3/4)}] {
M error\\
1 1.000e+00\\
2 7.322e-01\\
3 7.274e-01\\
5 6.287e-01\\
8 5.158e-01\\
13 1.685e-01\\
21 1.585e-01\\
34 5.964e-02\\
55 5.314e-02\\
89 2.388e-02\\
144 2.158e-02\\
233 1.847e-02\\
377 1.790e-02\\
610 7.939e-03\\
987 7.659e-03\\
1597 4.402e-03\\
2584 4.273e-03\\
4181 2.040e-03\\
6765 1.974e-03\\
10946 9.793e-04\\
17711 9.469e-04\\
28657 4.636e-04\\
46368 4.474e-04\\
75025 2.182e-04\\
121393 2.100e-04\\
196418 1.087e-04\\
317811 1.048e-04\\
514229 5.211e-05\\
832040 5.021e-05\\
1346269 2.541e-05\\
2178309 2.449e-05\\
3524578 1.235e-05\\
5702887 1.190e-05\\
9227465 6.009e-06\\
14930352 5.792e-06\\
24157817 2.925e-06\\
39088169 2.820e-06\\
63245986 1.419e-06\\
102334155 1.368e-06\\
165580141 6.900e-07\\
267914296 6.652e-07\\
433494437 3.352e-07\\
};
\addlegendentry{Fibonacci $\ell_\infty$-ball}
\addplot[brown,mark=diamond*,mark size=2.5] table[row sep=\\,x=M, y expr = { \thisrow{error} * (\thisrow{M})^(3/4) }] {
M error\\
1 1.000e+00\\
3 7.274e-01\\
8 5.148e-01\\
28 1.181e-01\\
104 2.457e-02\\
400 1.528e-02\\
1568 5.131e-03\\
6208 1.584e-03\\
24704 5.992e-04\\
98560 2.079e-04\\
393728 7.520e-05\\
1573888 2.602e-05\\
6293504 9.097e-06\\
25169920 3.286e-06\\
100671488 1.149e-06\\
402669568 4.070e-07\\
};
\addlegendentry{Korobov $\ell_\infty$-ball}
\end{loglogaxis}
\end{tikzpicture}
\caption{Scaled $L_2(\T^2)$ sampling error and number of sampling points for the approximation of the kink function $g$ from \eqref{equ:kink},
where $\mathrm{err}:=\|g-S^{\Lambda(\zb z,M)}_{I^2_N}g|L_2(\T^2)\|$.}
\label{figure:kink:linf_scaled}
\end{figure}

\subsection*{Acknowledgements}
\begin{sloppypar}
The authors acknowledge the fruitful discussions with A.~Hinrichs, M.~Ullrich and R.~Bergmann on this topic, especially
at the conference ``Approximationsmethoden und schnelle Algorithmen" in Hasenwinkel,
2014. Furthermore the authors thank V.N.~Temlyakov for his valuable comments and historical hints on that topic. The authors especially thank A. Hinrichs for pointing out an alternative proof argument for the non-optimality of rank-1 lattice sampling, cf. Remark \ref{rem:minkowski}. Moreover, LK and TV gratefully acknowledge support by the German Research Foundation (DFG) within the Priority Program
1324, project PO~711/10-2. Additionally, GB and TU acknowledge the support by the DFG Emmy-Noether programme (UL403/1-1)
and the Hausdorff-Center for Mathematics, University of Bonn.
\end{sloppypar}


\end{document}